\documentclass[11pt]{amsart}
\usepackage{amsmath,amssymb,amsthm}
\usepackage[margin=1in]{geometry}
\usepackage{enumitem}
\usepackage{booktabs}
\usepackage{array}
\usepackage[dvipsnames]{xcolor}
\usepackage[colorlinks=true,linkcolor=blue!50!black,citecolor=blue!50!black,urlcolor=blue!50!black]{hyperref}

\usepackage{hyperref}
\hypersetup{colorlinks=true,linkcolor=blue,citecolor=blue,urlcolor=blue}

\theoremstyle{plain}
\newtheorem{theorem}{Theorem}[section]
\newtheorem{lemma}[theorem]{Lemma}
\newtheorem{proposition}[theorem]{Proposition}
\newtheorem{corollary}[theorem]{Corollary}

\newtheorem{fact}[theorem]{Fact}

\theoremstyle{definition}
\newtheorem{definition}[theorem]{Definition}
\newtheorem{example}[theorem]{Example}

\newtheorem{problem}[theorem]{Problem}

\newcommand{\Z}{\mathbb{Z}}
\newcommand{\G}{\Gamma}
\newcommand{\Hq}{\mathcal{H}}
\newcommand{\Ev}{\mathcal{E}}
\newcommand{\Od}{\mathcal{O}}
\newcommand{\Bc}{\mathfrak{B}}
\newcommand{\gen}{\operatorname{gen}}
\newcommand{\genus}{\operatorname{g}}
\newcommand{\diam}{\operatorname{diam}}
\newcommand{\sdim}{\operatorname{sdim}}

\newcommand{\Q}{\mathrm{Q}}
\newcommand{\ord}[1]{o(#1)}
\newcommand{\cl}[1]{[#1]}

\title{On the prime index elements graph of finite groups}

\author{Pallabi Manna}
\address{Pallabi Manna\\ Department of Mathematics, Indian Institute of Technology Hyderabad, Kandi, Sangareddy 502284, Telangana, India}
\email{mannapallabimath001@gmail.com}

\author{Santanu Mandal}
\address{Santanu Mandal\\Department of Computer Science and Engineering, SRM University-AP, Amaravati, Andhra Pradesh 522240, India}
\email{santanu.vumath@gmail.com}


\begin{document}
	\maketitle

	\begin{abstract}
		For a finite group $G$, the prime index elements graph $\G(G)$, introduced by Anitha and Rajkumar, has vertex set $G$, two distinct elements being adjacent if the cyclic subgroup generated by one has prime index in that generated by the other. We show that $\G(G)$ is a bipartite blow-up of a weighted graph $\Hq(G)$ on the cyclic subgroups of $G$, which is a cartesian product of trees when $G$ is nilpotent, and use this to study several graph invariants. We characterize the groups for which $\G(G)$ is a star, a cograph, a chordal graph, or has a perfect code. We compute the sharp bounds for the independence and domination number. We give an exact formula for the strong metric dimension, bound the metric and strong metric dimensions in terms of the quotient graph, show that the upper bounds for the strong metric dimension are attained for cyclic p-groups, and all upper bounds for groups whose prime index elements graph is complete bipartite, and determine the groups attaining the lower bounds: the strong metric dimension attains it only for $\Z_p$ with $p$ odd, and we give a criterion for the metric dimension, which we solve completely for cyclic groups. For nilpotent $G$ we give necessary conditions for $\G(G)$ to be Hamiltonian, in terms of the numbers of covers of cyclic subgroups in the Sylow subgroups, and prove that they are also sufficient when all Sylow subgroups but one are cyclic. We also show that a nilpotent group with a Hamiltonian path is $\Z_3$ or of the form $K\times\Z_2$ with $|K|$ odd, and that in the abelian case above a Hamiltonian path exists exactly when a Hamiltonian cycle does. Finally, we characterize the finite groups $G$ for which $\G(G)$ is planar, toroidal or of genus two.
	\end{abstract}

	\section{Introduction}
	
	Associating graphs with groups is a classical way of encoding group-theoretic information in combinatorial form. The commuting graph appeared already in the work of Brauer and Fowler~\cite{BF} on groups of even order, and the prime graph (Gruenberg--Kegel graph) of a finite group played an important role in the study of finite simple groups~\cite{Williams}. Kelarev and Quinn~\cite{KQ} introduced the directed power graph, and Chakrabarty, Ghosh and Sen~\cite{CGS} its undirected version; see~\cite{KSCT} for a survey. Further examples include the non-commuting graph~\cite{AAM}, the enhanced power graph~\cite{AACNS}, the reduced power graph~\cite{RA17}, and the prime index graph of Akbari et al.~\cite{AAHRS}, whose vertices are the subgroups of $G$, two subgroups being adjacent when one has prime index in the other. Cameron's survey~\cite{Cameron} gives a unified account of many of these graphs. Research in this area has two aims: to recover properties of a group from its graph, and to produce interesting families of graphs from groups.
	
	Anitha and Rajkumar~\cite{AR} introduced an element-level analogue of the prime index graph.
	
	\begin{definition}\label{def:graph}
		Let $G$ be a finite group. The \emph{(undirected) prime index elements graph} $\G(G)$ is the simple graph with vertex set $G$ in which distinct $x,y\in G$ are adjacent if and only if $\langle y\rangle\subsetneq\langle x\rangle$ and $[\langle x\rangle:\langle y\rangle]$ is prime, or $\langle x\rangle\subsetneq\langle y\rangle$ and $[\langle y\rangle:\langle x\rangle]$ is prime.
	\end{definition}
	
	The graph $\G(G)$ is a spanning subgraph of the reduced power graph and hence of the power graph. It keeps only the edges between $x$ and $y$ for which $\langle x\rangle$ and $\langle y\rangle$ are adjacent in the poset of cyclic subgroups. Thus $\G(G)$ records the Hasse diagram of the poset of cyclic subgroups together with the number of generators of each cyclic subgroup, and it is much sparser than the power graph. Anitha and Rajkumar~\cite{AR} studied its basic properties, such as connectedness, bipartiteness and diameter, characterized the groups for which it is complete bipartite, and related it to the power graph and the reduced power graph.
	
	The importance of the present paper lies in a systematic treatment of $\G(G)$ based on a single structural observation (Section~\ref{sec:prelim}). Adjacency depends only on the cyclic subgroups generated, so $\G(G)$ is obtained from a weighted graph $\Hq(G)$ on the cyclic subgroups by replacing each cyclic subgroup $C$ by an independent set of its $\varphi(|C|)$ generators. For nilpotent $G$, $\Hq(G)$ is a Cartesian product of trees. This turns questions about $\G(G)$ into weighted problems on a much smaller graph. It lets us
	\begin{itemize}[leftmargin=*]
		\item characterize stars, cographs, chordal graphs and graphs with perfect codes;
		\item obtain sharp bounds for the independence and domination numbers, and determine when they are attained in several cases;
		\item give an exact formula for the strong metric dimension and bounds for the metric and strong metric dimensions in terms of the quotient graph $\Hq(G)$, show that the upper bounds are attained for cyclic $p$-groups and for all groups whose graph is complete bipartite, and characterize the groups attaining the lower bounds (completely for cyclic groups);
		\item give necessary conditions for Hamiltonicity of $\G(G)$ for nilpotent $G$, show that they are sufficient when all Sylow subgroups but one are cyclic (with an explicit numerical criterion when the remaining one is abelian), and treat Hamiltonian paths;
		\item characterize the groups for which $\G(G)$ is planar, toroidal or of genus two.
	\end{itemize}
	The genus results are the deepest part of the paper. They combine a reduction to a single component of an auxiliary ``commuting'' graph between involutions and subgroups of order $3$ with an orientation-balance argument for embeddings of bipartite graphs, and they give a complete topological classification of the kind available for power graphs and reduced power graphs. 
	
	\paragraph{Reading the group from the graph.} Many properties of $\G(G)$ translate into transparent conditions on $G$, so the graph can be used to recognize structural features of the group.
	\begin{itemize}[leftmargin=*]
		\item \emph{Element orders.} The neighbourhood of $e$ is the set of elements of prime order, and the bipartition of $\G(G)$ is given by the parity of the number of prime factors of $\ord x$ (Lemma~\ref{lem:bip}). The graph is a star, has a perfect code, or has independence number $|G|-1$ exactly when every non-identity element has prime order or $G\cong\Z_4,\Q_8$ (Theorems~\ref{thm:star},~\ref{thm:pc} and~\ref{thm:indep}). It is a tree exactly when every non-identity element has order $4$ or a prime order (Theorem~\ref{thm:chordal}).
		\item \emph{Cyclic subgroups.} The twin classes of generators record the cyclic subgroups and their numbers of generators. The metric dimension is at least $|G|-c(G)$, where $c(G)$ is the number of cyclic subgroups, and the strong metric dimension attains this value only for $\Z_p$ (Theorems~\ref{thm:dim} and~\ref{thm:lowattain}); for cyclic groups, attaining the lower bound is decided by the exponents in the prime factorization of $|G|$ (Theorem~\ref{thm:lowcyc}).
		\item \emph{Sylow structure.} For nilpotent $G$, a Hamiltonian cycle forces the Sylow $2$-subgroup to be $\Z_2$ and bounds the number of subgroups of order $p$ of each Sylow $p$-subgroup in terms of the order of its complement (Theorem~\ref{thm:hamnec}). For $G=P\times R\times\Z_2$ with $P$ an abelian $p$-group and $R$ cyclic, it is equivalent to an inequality between the rank of $P$ and $|R|$ (Corollary~\ref{cor:hamab}).
		\item \emph{Centralizers.} $\G(G)$ is planar if and only if every element order lies in $\{1,2,3,4,8,9\}\cup\{p\ge5\}$, and toroidal if and only if, in addition to the orders lying in $\Pi$, $G$ has exactly one cyclic subgroup of order $6$. Genus two is equivalent to explicit conditions on the centralizers of involutions and of subgroups of order $3$ (Theorems~\ref{thm:planar},~\ref{thm:torus} and~\ref{thm:genus2}).
	\end{itemize}
	Thus element orders, the poset of cyclic subgroups, the Sylow structure and centralizer data can each be detected by a graph invariant of $\G(G)$.
	
	\paragraph{Significance.} Among the graphs on a group defined through cyclic subgroups, $\G(G)$ is the sparsest: it is a spanning subgraph of the reduced power graph, which in turn is a spanning subgraph of the power graph. Because of this sparsity, its invariants respond to finer group-theoretic data than those of the denser graphs, such as the number of covers of a cyclic subgroup or the centralizer of a single involution. On the other hand, the blow-up description keeps every question finite and explicit. The method of this paper is to reduce a graph invariant to a weighted problem on the quotient $\Hq(G)$ and then to group theory. It applies equally to other graphs in which adjacency depends only on the cyclic subgroups generated, such as power graphs and enhanced power graphs, which are blow-ups by cliques instead of independent sets. To our knowledge, the metric dimensions, domination number, Hamiltonicity and genus of $\G(G)$ have not been studied yet. The results also leave natural open problems (Problem~\ref{prob:ham}).
	
	The paper is organized as follows. Section~\ref{sec:prelim} establishes the blow-up structure. Sections~\ref{sec:classes}--\ref{sec:genus} treat, in turn, special graph classes, the independence number, the domination number, the metric dimensions, Hamiltonicity and the genus. Throughout, all groups are finite, $D_n$ denotes the dihedral group of order $2n$, and $\Q_{2^m}$ the generalized quaternion group of order $2^m$, $m\ge3$.
	\section{Preliminaries and basic structure}\label{sec:prelim}
	
	\subsection{Notation}
	All groups are finite. For $x\in G$, $\ord{x}$ is its order. For $n\ge1$, $\Omega(n)$ denotes the number of prime factors of $n$ counted with multiplicity. For a graph $X$: $N(v)$, $N[v]=N(v)\cup\{v\}$, $\deg v$, $\delta(X)$, $\alpha(X)$, $\beta(X)$ (vertex cover number), $\nu(X)$ (matching number), $\gamma(X)$ (domination number), $\diam(X)$ and $\genus(X)$ (orientable genus) have their usual meanings. A vertex is \emph{universal} if it is adjacent to all other vertices.
	
	\subsection{Group-theoretic facts}
	We use the following standard results.
	\begin{fact}\label{fact:group}
		\begin{enumerate}[label=(\roman*)]
			\item (Cauchy) If a prime $p$ divides $|G|$, then $G$ has an element of order $p$.
			\item (Frobenius; see~\cite[Thm.~9.1.1]{Hall}) If $n$ divides $|G|$, the number of solutions of $x^n=1$ in $G$ is a multiple of $n$. Consequently a group of even order has an odd number of involutions, and if $3\mid|G|$ the number of elements of order $3$ is $\equiv 2\pmod 3$.
			\item (see~\cite[Thm.~5.4.10(ii)]{Gor}) A finite $p$-group with a unique subgroup of order $p$ is cyclic, or $p=2$ and it is generalized quaternion. The group $\Q_{2^m}$ has exponent $2^{m-1}$.
			\item The number of subgroups of order $p$ in a nontrivial finite $p$-group is $\equiv 1\pmod p$.
			\item A subgroup $S$ of a nilpotent group $G=P_1\times\cdots\times P_k$ ($P_i$ the Sylow subgroups) satisfies $S=\prod_i(S\cap P_i)$.
		\end{enumerate}
	\end{fact}
	
	\subsection{Graph-theoretic facts}
	\begin{fact}\label{fact:graph}
		\begin{enumerate}[label=(\roman*)]
			\item (Gallai) $\alpha(X)+\beta(X)=|V(X)|$. (K\H{o}nig) If $X$ is bipartite then $\beta(X)=\nu(X)$.
			\item (Youngs~\cite{Youngs}) A minimum-genus embedding of a connected graph is cellular, and a cellular embedding in the orientable surface $S_g$ satisfies $V-E+F=2-2g$.
			\item (Battle--Harary--Kodama--Youngs~\cite{BHKY}) The genus of a graph is the sum of the genera of its blocks.
			\item (Ringel~\cite{Ringel}) $\genus(K_{m,n})=\lceil (m-2)(n-2)/4\rceil$.
			\item If $X$ is connected and triangle-free then $\genus(X)\ge \lceil (E-2V+4)/4\rceil$.
		\end{enumerate}
	\end{fact}
	
	\subsection{The quotient graph}
	For a cyclic subgroup $C\le G$ let $\gen(C)$ be its set of generators; $|\gen(C)|=\varphi(|C|)$. The sets $\gen(C)$ partition $G$; we call them \emph{classes}, and write $\cl{x}=\gen(\langle x\rangle)$. Let $c(G)$ be the number of cyclic subgroups and $m(G)$ the number of maximal cyclic subgroups of $G$.
	
	\begin{definition}
		$\Hq(G)$ is the graph whose vertices are the cyclic subgroups of $G$, with $C\sim D$ if and only if one of $C,D$ is a subgroup of prime index in the other. Each vertex $C$ carries the weight $w(C)=\varphi(|C|)$.
	\end{definition}
	
	\begin{lemma}[Blow-up]\label{lem:blowup}
		\begin{enumerate}[label=(\roman*)]
			\item Two distinct elements of one class are non-adjacent and have the same neighbourhood.
			\item For distinct cyclic subgroups $C,D$, either every element of $\gen(C)$ is adjacent to every element of $\gen(D)$ (if and only if $C\sim D$ in $\Hq(G)$), or no such pair is adjacent.
		\end{enumerate}
		Thus $\G(G)$ is obtained from $\Hq(G)$ by replacing each vertex $C$ by an independent set of size $w(C)$ and each edge by a complete bipartite join. Moreover, for distinct $x,x'$ in one class, the transposition $(x\ x')$ is an automorphism of $\G(G)$.
	\end{lemma}
	\begin{proof}
		The adjacency condition involves only $\langle x\rangle$ and $\langle y\rangle$; if $\langle x\rangle=\langle y\rangle$ the index is $1$, which is not prime. The last claim follows from (i).
	\end{proof}
	
	\begin{lemma}\label{lem:conn}
		$\G(G)$ is connected.
	\end{lemma}
	\begin{proof}
		For $x\ne e$ choose a maximal subgroup $M$ of the cyclic group $\langle x\rangle$; its index is prime, so a generator $y$ of $M$ satisfies $y\sim x$ and $\ord y<\ord x$. Iterating gives a path from $x$ to $e$.
	\end{proof}
	
	\begin{lemma}[Bipartiteness]\label{lem:bip}
		Let $\Ev(G)=\{x:\Omega(\ord x)\text{ even}\}$ and $\Od(G)=\{x:\Omega(\ord x)\text{ odd}\}$. Every edge of $\G(G)$ joins $\Ev(G)$ to $\Od(G)$; in particular $\G(G)$ is bipartite and triangle-free.
	\end{lemma}
	\begin{proof}
		If $\langle y\rangle$ has prime index $r$ in $\langle x\rangle$ then $\ord x=r\,\ord y$, so $\Omega(\ord x)=\Omega(\ord y)+1$.
	\end{proof}
	
	\begin{lemma}\label{lem:prime}
		Let $y$ have prime order. For $x\ne y$, $x\sim y$ if and only if $x=e$, or $y\in\langle x\rangle$ and $\ord x/\ord y$ is prime. In particular $N(e)$ is the set of elements of prime order, and no two elements of prime order are adjacent.
	\end{lemma}
	\begin{proof}
		The only proper subgroup of $\langle y\rangle$ is $1$, of prime index $\ord y$. For two elements of prime orders $p,q$, containment of $\langle x\rangle$ in $\langle y\rangle$ with prime index is impossible since $p/q$ is never a prime.
	\end{proof}
	
	\begin{definition}
		A \emph{down-neighbour} of $x$ is a neighbour $y$ with $\langle y\rangle$ of prime index in $\langle x\rangle$.
	\end{definition}
	
	\begin{lemma}\label{lem:down}
		The number of down-neighbours of $x$ is $d(x)=\sum_{r\mid \ord x,\ r\text{ prime}}\varphi(\ord x/r)$, and every edge is a down-edge for exactly one of its endpoints. Hence $|E(\G(G))|=\sum_{x\in G}d(x)$.
	\end{lemma}
	\begin{proof}
		The cyclic group $\langle x\rangle$ has exactly one maximal subgroup of order $\ord x/r$ for each prime $r\mid\ord x$, with $\varphi(\ord x/r)$ generators. For an edge $xy$ exactly one of $\langle x\rangle,\langle y\rangle$ is the larger.
	\end{proof}
	
	\begin{lemma}[Distances]\label{lem:dist}
		Let $|G|\ge2$, $x\in C$, $y\in D$. If $C\ne D$ then $d_{\G}(x,y)=d_{\Hq}(C,D)$; if $C=D$ and $x\ne y$ then $d_\G(x,y)=2$. In particular $\Hq(G)$ is connected.
	\end{lemma}
	\begin{proof}
		A path in $\G$ projects, by Lemma~\ref{lem:blowup}, to a walk of the same length in $\Hq$, and an $\Hq$-path lifts to a $\G$-path by choosing any element in each class. Twins have a common neighbour by Lemma~\ref{lem:conn}.
	\end{proof}
	
	\begin{lemma}\label{lem:tree}
		If $P$ is a finite $p$-group then $\Hq(P)$ is a tree rooted at $1$ whose weights sum to $|P|$; the neighbours of the root are the subgroups of order $p$.
	\end{lemma}
	\begin{proof}
		A cyclic subgroup $C\ne1$ of order $p^a$ has a unique maximal subgroup, its \emph{parent}. Every edge of $\Hq(P)$ joins a vertex to its parent, so $|E|=|V|-1$, and following parents reaches $1$, so $\Hq(P)$ is connected; hence it is a tree. The weights sum to $|P|$ since the classes partition $P$.
	\end{proof}
	
	\begin{lemma}\label{lem:product}
		If $G=P_1\times\cdots\times P_k$ is nilpotent with Sylow subgroups $P_i$, then
		\[ \Hq(G)\cong \Hq(P_1)\,\square\cdots\square\,\Hq(P_k)\quad\text{(Cartesian product)},\qquad w(C_1\times\cdots\times C_k)=\prod_i w(C_i). \]
	\end{lemma}
	\begin{proof}
		By Fact~\ref{fact:group}(v) every subgroup is $\prod(S\cap P_i)$, and it is cyclic if and only if each factor is (the factors have coprime orders). $\prod C'_i$ has prime index $p_i$ in $\prod C_i$ if and only if $C'_j=C_j$ for $j\neq i$ and $C'_i$ is maximal in $C_i$; this is adjacency in the Cartesian product. Finally $\varphi$ is multiplicative on coprime arguments.
	\end{proof}
	
	\section{Stars, cographs, chordal graphs and perfect codes}\label{sec:classes}
	
	\begin{lemma}\label{lem:star}
		Let $X$ be a connected triangle-free graph with at least one edge. Then $X$ has a universal vertex if and only if $X$ is a star if and only if  $\nu(X)=1$.
	\end{lemma}
	\begin{proof}
		If $u$ is universal, two other vertices cannot be adjacent (triangle), so $X$ is a star. A star has $\nu=1$. If $\nu(X)=1$, fix an edge $uv$; if there are edges $uw$ ($w\ne v$) and $vz$ ($z\ne u$), they must meet, forcing $w=z$ and a triangle. So all edges contain $u$, or all contain $v$; by connectedness $X$ is a star.
	\end{proof}
	
	\begin{theorem}\label{thm:star}
		Let $|G|\ge2$. The following are equivalent:
		\begin{enumerate}[label=(\alph*)]
			\item $\G(G)$ has a universal vertex;
			\item $\G(G)$ is a star $K_{1,|G|-1}$;
			\item $\nu(\G(G))=1$;
			\item every non-identity element of $G$ has prime order, or $G\cong\Z_4$, or $G\cong\Q_8$.
		\end{enumerate}
		In case (d) with every non-identity element of prime order and $|G|\ge3$, the unique universal vertex is $e$; for $\Z_4$ and $\Q_8$ it is the unique involution; for $\Z_2$ both vertices are universal.
	\end{theorem}
	\begin{proof}
		(a)$\Leftrightarrow$(b)$\Leftrightarrow$(c) is Lemma~\ref{lem:star} with Lemmas~\ref{lem:conn},~\ref{lem:bip}.
		
		(b)$\Rightarrow$(d). Let $\G(G)$ be a star with centre $u$. If $u=e$, every $x\ne e$ is adjacent to $e$, hence of prime order by Lemma~\ref{lem:prime}. If $u\ne e$, then $e$ is a leaf, so $u$ is the only element of prime order; as $u^{-1}$ also has prime order, $\ord u=2$. By Cauchy, $G$ is a $2$-group. For $y\notin\{e,u\}$, a generator of a maximal subgroup of $\langle y\rangle$ is a neighbour of $y$, hence equals $u$; so $\ord y=4$. Thus $G$ is a $2$-group of exponent $4$ with a unique involution, so $G\cong\Z_4$ or $\Q_8$ by Fact~\ref{fact:group}(iii).
		
		(d)$\Rightarrow$(b). If all non-identity elements have prime order, $e$ is adjacent to all of them and there are no other edges (Lemma~\ref{lem:prime}). In $\Z_4$ the involution $u$ is adjacent to $e$ and to both generators, and there are no other edges, so $\G(\Z_4)=K_{1,3}$. In $\Q_8$, $e\sim-1$; each of the six elements of order $4$ is adjacent to $-1$ (index $2$) and to nothing else, since distinct cyclic subgroups of order $4$ are incomparable; so $\G(\Q_8)=K_{1,7}$.
		
		The final statement is clear since a star $K_{1,n}$ with $n\ge2$ has a unique universal vertex.
	\end{proof}
	
	\begin{lemma}\label{lem:P4}
		A connected bipartite graph is $P_4$-free if and only if it is complete bipartite.
	\end{lemma}
	\begin{proof}
		In an induced $a$--$b$--$c$--$d$ the endpoints lie in different parts and are non-adjacent, which is impossible in a complete bipartite graph. Conversely, if $u,v$ in different parts are non-adjacent, a shortest $u$--$v$ path is induced, of odd length $\ge3$, and its first four vertices induce $P_4$.
	\end{proof}
	
	\begin{theorem}[Cographs]\label{thm:cograph}
		$\G(G)$ is a cograph if and only if it is complete bipartite. Consequently, $\G(G)$ is a cograph if and only if one of the following holds: (1) every non-identity element of $G$ has prime order; (2) $G\cong\Z_{p^2}$; (3) $G\cong\Z_{pq}$ for distinct primes $p,q$; (4) $G\cong\Q_8$. In these cases $\G(G)$ is respectively $K_{1,|G|-1}$, $K_{p-1,\,p^2-p+1}$, $K_{p+q-2,\,(p-1)(q-1)+1}$ and $K_{1,7}$.
	\end{theorem}
	\begin{proof}
		By Lemmas~\ref{lem:conn} and~\ref{lem:bip}, $\G(G)$ is connected and bipartite, so by Lemma~\ref{lem:P4} it is a cograph if and only if it is complete bipartite. The finite groups $G$ for which $\G(G)$ is complete bipartite are determined by Anitha and Rajkumar~\cite{AR}; they are exactly the groups in (1)--(4).
		The parts of $K_{s,t}$ are $\Ev(G)$ and $\Od(G)$ (Lemma~\ref{lem:bip}), which gives the stated values of $s$ and $t$: in case (1), $\Ev(G)=\{e\}$; for $\Z_{p^2}$, $\Od(G)$ consists of the $p-1$ elements of order $p$; for $\Z_{pq}$, of the $p+q-2$ elements of prime order; and for $\Q_8$, $\Od(G)=\{-1\}$.
	\end{proof}
	
	\begin{theorem}[Chordality]\label{thm:chordal}
		The following are equivalent: (a) $\G(G)$ is chordal; (b) $\G(G)$ is a tree; (c) every non-identity element of $G$ has order $4$ or a prime order. In that case $e$ is adjacent to every element of prime order, each element of order $4$ is a leaf attached to its square, and $\diam\G(G)\le4$.
	\end{theorem}
	\begin{proof}
		(a)$\Rightarrow$(b): in a chordal graph a shortest cycle is a triangle; $\G(G)$ is triangle-free and connected, hence a tree. (b)$\Rightarrow$(a) is trivial.
		
		(c)$\Rightarrow$(b): By Lemma~\ref{lem:down}, $d(x)=1$ for every $x\neq e$ (if $\ord x$ is prime the only down-neighbour is $e$; if $\ord x=4$ it is $x^2$), so $|E|=|G|-1$ and $\G(G)$ is a tree by Lemma~\ref{lem:conn}.
		
		(b)$\Rightarrow$(c): Let $\ord x=n\ge3$. Then $\varphi(n)\ge2$, so $x$ has a twin $x'$ with $N(x)=N(x')$; two distinct neighbours $y,z$ of $x$ would give the $4$-cycle $x\,y\,x'\,z$. Hence $\deg x=1$, so $d(x)=\sum_{r\mid n}\varphi(n/r)=1$, which forces $n$ to have a single prime divisor $r$ with $\varphi(n/r)=1$, i.e.\ $n=r$ or $n=4$.
		
		Elements of order $4$ have no up-neighbours (their order would be $8$ or $4r$), whence the final description.
	\end{proof}
	
	A perfect code $C\subseteq V$ in graph theory is a subset of vertices where every vertex in the graph is at a distance of at most one from exactly one vertex in that subset. The set $C$ is an independent set. A perfect code is also known as an efficient dominating set because every vertex outside $C$ is adjacent to precisely one vertex inside $C$.
	
	\begin{theorem}[Perfect codes]\label{thm:pc}
		$\G(G)$ has a perfect code if and only if every non-identity element of $G$ has prime order, or $G\cong\Z_4$ or $\Q_8$. The perfect codes are exactly the singletons $\{u\}$ with $u$ universal.
	\end{theorem}
	\begin{proof}
		If $|G|=1$ the claim is trivial; let $|G|\ge2$ and let $C$ be a perfect code. If $c\in C$ has $\ord c\ge3$, a twin $c'\ne c$ would have $N[c']=\{c'\}\cup N(c)$; $c'\in C$ would give $N[c]\cap N[c']\supseteq N(c)\ne\emptyset$, and $u\in N(c)\cap C$ would give $N[u]\cap N[c]\ni u$---both impossible, so $c'$ is not dominated exactly once. Hence $C\subseteq\{e\}\cup\{\text{involutions}\}$. Any two elements of this set have intersecting closed neighbourhoods ($e$ is adjacent to every involution, and two involutions share the neighbour $e$), so $|C|=1$ and $C=\{u\}$ with $N[u]=G$. Conversely a universal vertex is a perfect code. Now apply Theorem~\ref{thm:star}.
	\end{proof}
	
	\begin{corollary}\label{cor:families1}
		For a finite nilpotent group $G$:
		\begin{enumerate}[label=(\roman*)]
			\item $\G(G)$ has a universal vertex (equivalently, a perfect code) if and only if $G$ is a $p$-group of exponent $p$, or $G\cong\Z_4,\Q_8$;
			\item $\G(G)$ is a cograph if and only if $G$ is a $p$-group of exponent $p$, or $G\cong\Z_{p^2},\Z_{pq},\Q_8$;
			\item $\G(G)$ is chordal if and only if $G$ is a $p$-group of exponent $p$ with $p$ odd, or a $2$-group of exponent dividing $4$.
		\end{enumerate}
		For the symmetric, alternating and dihedral groups:
		\begin{center}
			\begin{tabular}{@{}llll@{}}
				\toprule
				& universal vertex / perfect code & cograph & chordal (tree)\\
				\midrule
				$S_n$ & $n\le3$ & $n\le3$ & $n\le4$\\
				$A_n$ & $n\le5$ & $n\le5$ & $n\le6$\\
				$D_n$ & $n=1$ or $n$ prime & $n=1$ or $n$ prime & $n\in\{1,2,4\}$ or $n$ an odd prime\\
				\bottomrule
			\end{tabular}
		\end{center}
	\end{corollary}
	\begin{proof}
		In a nilpotent group whose order has two distinct prime divisors $p,q$, commuting elements of orders $p$ and $q$ have product of order $pq$; this gives the nilpotent statements from Theorems~\ref{thm:star},~\ref{thm:cograph},~\ref{thm:chordal}. The element orders of $S_n$ for $n\le4$ are contained in $\{1,2,3,4\}$, while $S_n$ ($n\ge5$) contains $(1\,2)(3\,4\,5)$ of order $6$, and $S_4$ contains a $4$-cycle. The element orders of $A_4,A_5,A_6$ are $\{1,2,3\},\{1,2,3,5\},\{1,2,3,4,5\}$ (in $A_6$ an element of order $6$ would be a $6$-cycle or a $3$-cycle times a transposition, both odd), while $A_n$ ($n\ge7$) contains $(1\,2\,3)(4\,5)(6\,7)$ of order $6$. In $D_n$ reflections have order $2$ and the rotation subgroup is $\Z_n$. None of $S_n,A_n,D_n$ is isomorphic to $\Z_4$, $\Q_8$, $\Z_{p^2}$ or $\Z_{pq}$ except the cases already covered by condition (1) of the relevant theorem ($D_1\cong\Z_2$, $A_3\cong\Z_3$).
	\end{proof}
	
	\section{Independence number}\label{sec:indep}
	
	\begin{theorem}\label{thm:indep}
		Let $|G|\ge2$. Then
		\[
		\Big\lceil \tfrac{|G|}{2}\Big\rceil\ \le\ \max\{|\Ev(G)|,|\Od(G)|\}\ \le\ \alpha(\G(G))\ =\ |G|-\nu(\G(G))\ \le\ |G|-1 .
		\]
		Moreover:
		\begin{enumerate}[label=(\roman*)]
			\item $\alpha(\G(G))=|G|-1$ if and only if every non-identity element of $G$ has prime order, or $G\cong\Z_4,\Q_8$;
			\item $\alpha(\G(G))=|G|/2$ if and only if $\G(G)$ has a perfect matching; this holds for every $G=A\times\Z_2$ with $|A|$ odd.
		\end{enumerate}
	\end{theorem}
	\begin{proof}
		$\Ev,\Od$ partition $G$ and are independent (Lemma~\ref{lem:bip}); $\alpha=|G|-\nu$ by Fact~\ref{fact:graph}(i); and $\nu\ge1$ as $\G(G)$ is connected with at least two vertices. (i) follows from Theorem~\ref{thm:star}((c)$\Leftrightarrow$(d)). For (ii), $\alpha=|G|/2$ if and only if $\nu=|G|/2$. If $G=A\times\langle z\rangle$, $\ord z=2$, $|A|$ odd, and $a\in A$ with $m=\ord a$, then $\ord{az}=2m$ and $(az)^{m+1}=a$ since $m+1$ is even, so $\langle a\rangle$ has index $2$ in $\langle az\rangle$ and $a\sim az$. The edges $\{a,az\}$, $a\in A$, form a perfect matching.
	\end{proof}
	
	\begin{proposition}\label{prop:alphaH}
		$\alpha(\G(G))$ equals the maximum weight $\sum_{C\in\mathcal S}\varphi(|C|)$ of an independent set $\mathcal S$ of $\Hq(G)$.
	\end{proposition}
	\begin{proof}
		The classes meeting an independent set $I$ of $\G(G)$ are pairwise non-adjacent in $\Hq(G)$ (Lemma~\ref{lem:blowup}(ii)), so $|I|$ is at most the weight of this independent set; conversely a union of pairwise non-adjacent classes is independent.
	\end{proof}
	
	\section{Domination number}\label{sec:dom}
	
	For a set $\mathcal S$ of vertices of $\Hq(G)$, let $I(\mathcal S)$ be the set of $C\in\mathcal S$ having no neighbour in $\mathcal S$.
	
	\begin{theorem}[Reduction]\label{thm:domred}
		\[ \gamma(\G(G))=\min\Big\{|\mathcal S|+\sum_{C\in I(\mathcal S)}\big(w(C)-1\big)\ :\ \mathcal S \text{ a dominating set of } \Hq(G)\Big\}. \]
	\end{theorem}
	\begin{proof}
		Let $D$ dominate $\G(G)$ and let $\mathcal S$ be the set of classes meeting $D$. If $C\notin\mathcal S$, any $v\in\gen(C)$ has a neighbour in $D$, whose class is in $\mathcal S$ and adjacent to $C$; so $\mathcal S$ dominates $\Hq(G)$. If $C\in I(\mathcal S)$ and $v\in\gen(C)$, all neighbours of $v$ lie in classes outside $\mathcal S$, so $v\in D$; hence $\gen(C)\subseteq D$, and $|D|$ is at least the displayed quantity. Conversely, given $\mathcal S$, take one element from each class of $\mathcal S$ together with all elements of each class in $I(\mathcal S)$; Lemma~\ref{lem:blowup} shows this dominates $\G(G)$.
	\end{proof}
	
	\begin{theorem}\label{thm:dombounds}
		For $G\ne1$,
		\[ 1\le \gamma(\G(G))\le \min\{\,c(G)-m(G),\ |\Ev(G)|,\ |\Od(G)|\,\}, \]
		and $\gamma(\G(G))=1$ if and only if $G$ is as in Theorem~\ref{thm:star}(d). The bound $c(G)-m(G)$ is attained when every non-identity element has prime order, and by $\Z_{p^2}$ with $p$ odd.
	\end{theorem}
	\begin{proof}
		Choose a generator $x_C$ of every cyclic subgroup $C\ne1$ that is not maximal cyclic and set $D=\{e\}\cup\{x_C\}$, so $|D|=c(G)-m(G)$. Every class $C\ne1$ has a maximal subgroup $M$ whose generators are adjacent to all of $\gen(C)$; either $M=1$ and $e\in D$, or $M\ne1$ is properly contained in $C$, hence not maximal cyclic, and $x_M\in D$. Thus $D$ dominates. Since $\G(G)$ is connected and bipartite, each part dominates. $\gamma=1$ means a universal vertex exists. If all non-identity elements have prime order then $c-m=1$; for $\Z_{p^2}$ with $p$ odd, $c-m=2$, while $\gamma\ge2$ because $\G(\Z_{p^2})$ has no universal vertex (Theorem~\ref{thm:star}).
	\end{proof}
	
	\begin{example}
		For distinct primes $p,q$, $\gamma(\G(\Z_{pq}))=2<3=c-m$: if $\ord x=p$ then $\{e,x\}$ dominates (the generators of $\Z_{pq}$ are adjacent to $x$, all other elements to $e$), and there is no universal vertex.
	\end{example}
	
	\subsection{When is the bound $c(G)-m(G)$ attained?}
	Write $\mathcal N(G)$ for the set of cyclic subgroups of $G$ that are not maximal cyclic; thus $|\mathcal N(G)|=c(G)-m(G)$, and $1\in\mathcal N(G)$ when $G\ne1$. A \emph{cover} of a cyclic subgroup $C$ is a cyclic subgroup $D\supsetneq C$ with $[D:C]$ prime; so $C\in\mathcal N(G)$ if and only if $C$ has a cover. Recall that a cyclic group of prime power order has a unique maximal subgroup, and a cyclic group whose order has at least two prime divisors has at least two.
	
	\begin{theorem}\label{thm:domattain}
		Let $G\ne1$.
		\begin{enumerate}[label=(\roman*)]
			\item \emph{(Strict inequality.)} Suppose some $C\in\mathcal N(G)$, $C\ne1$, has no cover of prime power order. Then $\gamma(\G(G))\le c(G)-m(G)-1$. In particular this holds if some non-maximal cyclic subgroup of $G$ has order that is not a prime power.
			\item \emph{(Equality.)} Suppose that $G$ has a subgroup of prime order which is a maximal cyclic subgroup, and that every $C\in\mathcal N(G)\setminus\{1\}$ is the unique maximal subgroup of some maximal cyclic subgroup $M_C$ of $G$. Then $\gamma(\G(G))=c(G)-m(G)$.
			\item \emph{(Cyclic groups.)} $\gamma(\G(\Z_n))=c(\Z_n)-m(\Z_n)$ if and only if $n$ is prime or $n=p^2$ with $p$ an odd prime.
		\end{enumerate}
	\end{theorem}
	\begin{proof}
		(i) Choose a generator $x_{C'}$ of each $C'\in\mathcal N(G)\setminus\{1,C\}$, and put $D=\{e\}\cup\{x_{C'}\}$, so $|D|=c(G)-m(G)-1$. Let $v\notin D$ and $V=\langle v\rangle$. Then $V\ne1$, and each maximal subgroup $M$ of $V$ is either $1$ or a proper subgroup of the cyclic group $V$, hence in $\mathcal N(G)$; all elements of $\gen(V)$ are adjacent to all elements of $\gen(M)$. If $V$ has a maximal subgroup $M\ne C$, then $e$ or $x_M$ lies in $D$ and is adjacent to $v$. Otherwise $C$ is the unique maximal subgroup of $V$, so $V$ is a cover of $C$ of prime power order, which is excluded. (Note that $V=C$ is covered by the first case, since $C\ne1$ has a maximal subgroup different from $C$.) Hence $D$ dominates. If $|C|$ is not a prime power, then no cover of $C$ has prime power order.
		
		(ii) Let $U_0$ be a subgroup of prime order that is maximal cyclic. For each $T\in\{U_0\}\cup\{M_C: C\in\mathcal N(G)\setminus\{1\}\}$ the closed neighbourhood of $T$ in $\Hq(G)$ is $N_\Hq[U_0]=\{U_0,1\}$, respectively $N_\Hq[M_C]=\{M_C,C\}$, since a maximal cyclic subgroup has no cover and its only maximal subgroup is $1$, respectively $C$. These $c(G)-m(G)$ sets are pairwise disjoint: $M_C$ determines $C$; maximal cyclic subgroups do not lie in $\mathcal N(G)$; and $U_0\ne M_C$ because their maximal subgroups $1$ and $C$ differ. Pick one element $y_T$ in each such class $T$. By Lemma~\ref{lem:blowup}, $N_\G[y_T]$ is contained in the union of the classes in $N_\Hq[T]$; these unions are pairwise disjoint, and every dominating set meets each $N_\G[y_T]$. Hence $\gamma(\G(G))\ge c(G)-m(G)$, and equality follows from Theorem~\ref{thm:dombounds}.
		
		(iii) Let $n=p^k$. Then $c(\Z_n)=k+1$ and $m(\Z_n)=1$, so $c-m=k$. For $k=1$, $\gamma=1=k$. For $k=2$, Theorem~\ref{thm:star} shows that $\G(\Z_{p^2})$ has a universal vertex if and only if $p=2$; so $\gamma(\G(\Z_4))=1<2$, while $\gamma(\G(\Z_{p^2}))=2$ for odd $p$. For $k\ge3$, let $x_j$ be a generator of the subgroup $L_j$ of order $p^j$. Then $\{x_1,\dots,x_{k-1}\}$ dominates $\G(\Z_n)$: $e$ is adjacent to $x_1$; every other generator of $L_1$ is adjacent to $x_2$; and every element of order $p^j$ with $j\ge2$ is adjacent to $x_{j-1}$. Hence $\gamma\le k-1<c-m$. If $n$ has two distinct prime divisors $p,q$, let $C$ be the subgroup of order $n/p$. Then $C\ne1$, its only cover is $\Z_n$, and $n$ is not a prime power; so (i) gives $\gamma<c-m$.
	\end{proof}
	
	\section{Metric and strong metric dimension}\label{sec:metric}
	
	A set $W\subseteq V(X)$ \emph{resolves} a connected graph $X$ if the vectors $(d(v,w))_{w\in W}$ are pairwise distinct; $\dim(X)$ is the least size of a resolving set. A vertex $s$ \emph{strongly resolves} $u\ne v$ if $d(v,s)=d(v,u)+d(u,s)$ or $d(u,s)=d(u,v)+d(v,s)$; $\sdim(X)$ is the least size of a set strongly resolving all pairs. Write $u\lhd v$ ($u$ is \emph{maximally distant} from $v$) if $d(v,w)\le d(v,u)$ for all $w\in N(u)$; $u,v$ are MMD if $u\lhd v$ and $v\lhd u$. The \emph{strong resolving graph} $X_{SR}$ joins MMD pairs.
	
	\begin{lemma}\label{lem:nD}
		(cf.~\cite{CEJO}) For a connected graph $X$ with $n$ vertices and diameter $D$, $\dim(X)\le\sdim(X)\le n-D$.
	\end{lemma}
	\begin{proof}
		If $s$ strongly resolves $u\ne v$, say $d(v,s)=d(v,u)+d(u,s)$, then $d(v,s)>d(u,s)$, so $s$ resolves them. Let $v_0,\dots,v_D$ be a geodesic. The set $V\setminus\{v_1,\dots,v_D\}$ resolves ($v_0$ separates the $v_i$), and $V\setminus\{v_0,\dots,v_{D-1}\}$ strongly resolves (pairs meeting it are resolved by their own member; $v_j$ lies on the $v_i$--$v_D$ geodesic for $i<j<D$).
	\end{proof}
	
	\begin{lemma}[Oellermann--Peters-Fransen~\cite{OPF}]\label{lem:OPF}
		For a connected graph $X$, a set is strong resolving if and only if it is a vertex cover of $X_{SR}$; hence $\sdim(X)=|V(X)|-\alpha(X_{SR})$.
	\end{lemma}
	
	\begin{theorem}\label{thm:dim}
		For $|G|\ge2$: \ $|G|-c(G)\le\dim(\G(G))\le|G|-c(G)+\dim(\Hq(G))\le |G|-\diam(\Hq(G))$.
	\end{theorem}
	\begin{proof}
		If a resolving set missed two twins $u,v$, the twin automorphism (Lemma~\ref{lem:blowup}) shows $d(u,w)=d(v,w)$ for all $w\notin\{u,v\}$; so it contains all but at most one element of each class. For the upper bound let $R$ resolve $\Hq(G)$, choose a representative $x_C$ of each class, and let $W$ consist of all non-representatives and all $x_C$ with $C\in R$. Two representatives $x_C,x_D\notin W$ ($C\ne D$, $C,D\notin R$) are separated by $x_E$, where $E\in R$ resolves $C,D$ in $\Hq$ (Lemma~\ref{lem:dist}). The last inequality is Lemma~\ref{lem:nD} applied to $\Hq(G)$.
	\end{proof}
	
	\begin{lemma}\label{lem:MMD}
		(a) Distinct twins are MMD. (b) Let $u\in\gen(C)$, $v\in\gen(D)$, $C\ne D$. Then $u\lhd v$ if and only if (b1) $C\sim D$ in $\Hq$ implies $w(D)=1$, and (b2) every $\Hq$-neighbour $W\ne D$ of $C$ satisfies $d_\Hq(W,D)\le d_\Hq(C,D)$.
		(c) Whether elements of two classes are MMD does not depend on the chosen elements.
	\end{lemma}
	\begin{proof}
		(a) $d(u,v)=2$ and $N(u)=N(v)$. (b) By Lemma~\ref{lem:dist}, a neighbour $w$ of $u$ in a class $W\neq D$ has $d(v,w)=d_\Hq(W,D)$; if $C\sim D$ then $d(u,v)=1$ and the other elements of $\gen(D)$ are neighbours of $u$ at distance $2$ from $v$. (c) follows from the twin automorphisms.
	\end{proof}
	
	Let $Q(G)$ be the graph on the cyclic subgroups in which $C\ne D$ are adjacent if and only if their elements are MMD (well defined by Lemma~\ref{lem:MMD}(c)).
	
	\begin{theorem}\label{thm:sdim}
		$\sdim(\G(G))=|G|-\alpha(Q(G))$. Consequently
		\[ |G|-c(G)\le\dim(\G(G))\le\sdim(\G(G))\le\min\{\,|G|-c(G)+\sdim(\Hq(G)),\ |G|-\diam(\G(G))\,\}. \]
	\end{theorem}
	\begin{proof}
		By Lemma~\ref{lem:MMD}(a), each class of size $\ge2$ is a clique of $\G_{SR}$, so an independent set of $\G_{SR}$ contains at most one element per class, and these classes must be pairwise non-adjacent in $Q(G)$; hence $\alpha(\G_{SR})=\alpha(Q(G))$, and Lemma~\ref{lem:OPF} applies. By Lemma~\ref{lem:MMD}(b), MMD classes are MMD in $\Hq$, so $Q(G)$ is a subgraph of $\Hq_{SR}$ and $\alpha(Q)\ge\alpha(\Hq_{SR})=c(G)-\sdim(\Hq)$. The rest is Theorem~\ref{thm:dim} and Lemma~\ref{lem:nD}.
	\end{proof}
	
	\begin{theorem}\label{thm:dimcyc}
		Let $G=\Z_{p^n}$.
		\begin{enumerate}[label=(\roman*)]
			\item $\sdim(\G(G))=p^n-n$ if $n\ge2$ or $p=2$; and $\sdim(\G(\Z_p))=p-2$ for odd $p$.
			\item For odd $p$: $\dim(\G(G))=p^n-n$ if $n\in\{2,3\}$, and $p^n-n-1$ otherwise.
			\item For $p=2$: $\dim(\G(G))=2^n-n$ if $n\in\{1,2,3,4\}$, and $2^n-n-1$ if $n\ge5$.
		\end{enumerate}
	\end{theorem}
	\begin{proof}
		$\Hq$ is the path $L_0\cdots L_n$ with $|L_0|=1$ and $|L_i|=\varphi(p^i)$, so $|L_i|\ge2$ for $i\ge1$ unless $p=2$, $i=1$.
		
		(i) Let $i<j$, $u\in L_i$, $v\in L_j$. By Lemma~\ref{lem:MMD}(b): if $i\ge1$, the neighbour class $L_{i-1}$ violates (b2), so $u\lhd v$ forces $i=0$, and then $u\lhd v$ if and only if $j\ge2$ or $|L_1|=1$. Similarly $v\lhd u$ forces $j=n$, and then holds if and only if $i<n-1$ or $|L_i|=1$. Hence $Q(G)$ has the single edge $L_0L_n$ when $n\ge2$ (so $\alpha(Q)=n$), $Q=K_2$ when $n=1$, $p=2$ (so $\alpha=1$), and no edge when $n=1$, $p$ odd (so $\alpha=2$). Apply Theorem~\ref{thm:sdim}.
		
		(ii),(iii) A resolving set of size $|G|-c(G)=p^n-n-1$ omits exactly one element $x_i\in L_i$ for each $i$; all such sets are equivalent under twin automorphisms, so they all resolve or none does, and they resolve if and only if the $x_i$ are pairwise resolved. Let $J=\{j:|L_j|\ge2\}$ (so $J=\{1,\dots,n\}$ for odd $p$, $J=\{2,\dots,n\}$ for $p=2$). A landmark in $L_j$, $j\in J$, is at distance $f_i(j)=|i-j|$ from $x_i$ if $i\ne j$, and $2$ if $i=j$. For $i<k$: $f_i(j)=f_k(j)$ for $j\notin\{i,k\}$ if and only if $j=(i+k)/2$, while at $j\in\{i,k\}$ equality holds if and only if $k-i=2$. If $k-i\ne2$, the pair is separated by $j=k$ or $j=i$ whenever these lie in $J$; otherwise $p=2$, $(i,k)=(0,1)$, which is separated by $j=2$ when $n\ge2$ and is unresolved when $n=1$. If $k-i=2$, the pair is resolved if and only if some $j\in J$ satisfies $j<i$ or $j>i+2$. For odd $p$ this fails exactly for $(n,i)=(2,0),(3,1)$; for $p=2$ exactly for $n=i+2\le 4$. In the failing cases $\dim\ge p^n-n$, and $\dim\le\sdim=p^n-n$ by (i).
	\end{proof}
	
	\subsection{Groups attaining the upper bounds}
	We now show that the upper bounds of Theorems~\ref{thm:dim} and~\ref{thm:sdim} are attained, i.e.\ they cannot be improved in general. We use the following elementary fact: if $X$ is a connected graph with $n\ge3$ vertices that is not complete, then $\diam(X)\ge2$, and Lemma~\ref{lem:nD} gives
	\begin{equation}\label{eq:n-2}
		\dim(X)\ \le\ \sdim(X)\ \le\ n-2 .
	\end{equation}
	
	\begin{theorem}\label{thm:metricattain}
		\begin{enumerate}[label=(\roman*)]
			\item Let $G$ be a finite group with $c(G)\ge3$ such that $\G(G)$ is complete bipartite, i.e.\ (Theorem~\ref{thm:cograph}) every non-identity element of $G$ has prime order and $G$ is not cyclic of prime order, or $G\cong\Z_{p^2}$, $\Z_{pq}$ or $\Q_8$. Then
			\begin{align*}
				\dim(\G(G))=\sdim(\G(G))&=|G|-c(G)+\dim(\Hq(G))=|G|-c(G)+\sdim(\Hq(G))\\
				&=|G|-\diam(\G(G))=|G|-2 .
			\end{align*}
			So all the upper bounds of Theorems~\ref{thm:dim} and~\ref{thm:sdim} are attained. In particular this holds for $S_3$, $A_4$, $A_5$, for every $p$-group of exponent $p$ other than $\Z_p$, and for the dihedral groups $D_p$, $p$ prime.
			\item Let $G=\Z_{p^n}$ with $n\ge2$. Then
			\[ \sdim(\G(G))=|G|-c(G)+\sdim(\Hq(G))=|G|-\diam(\G(G))=p^n-n, \]
			so both upper bounds of Theorem~\ref{thm:sdim} are attained. Moreover $\dim(\G(G))=|G|-c(G)+\dim(\Hq(G))$ if and only if $p$ is odd and $n\in\{2,3\}$, or $p=2$ and $n\in\{2,3,4\}$; for all other $n\ge2$, $\dim(\G(G))=|G|-c(G)$ attains the lower bound of Theorem~\ref{thm:dim}.
		\end{enumerate}
	\end{theorem}
	\begin{proof}
		(i) Let $n=|G|$ and $\G(G)=K_{s,t}$ with $s+t=n$. Two vertices in the same part have the same neighbourhood (the other part) and are non-adjacent, so a vertex $w$ outside both of them has the same distance to each. Hence a resolving set contains all but at most one vertex of each part, and $\dim(\G(G))\ge n-2$. Since $\G(G)$ is connected, bipartite and $n\ge3$, it is not complete and has diameter $2$, so \eqref{eq:n-2} gives $\dim(\G(G))=\sdim(\G(G))=n-2=n-\diam(\G(G))$.
		
		The graph $\Hq(G)$ is connected (Lemma~\ref{lem:dist}) and bipartite (Lemma~\ref{lem:bip}, since all elements of a class have the same order), and it has $c=c(G)\ge3$ vertices; hence it is not complete, and \eqref{eq:n-2} gives $\dim(\Hq(G))\le\sdim(\Hq(G))\le c-2$. Combining with Theorems~\ref{thm:dim} and~\ref{thm:sdim},
		\[ n-2=\dim(\G(G))\le n-c+\dim(\Hq(G))\le n-c+\sdim(\Hq(G))\le n-2, \]
		and similarly $n-2=\sdim(\G(G))\le n-c+\sdim(\Hq(G))\le n-2$. So equality holds throughout. The listed groups satisfy the hypothesis: their non-identity elements have prime order, and $c(G)\ge3$ unless $G\cong\Z_p$.
		
		(ii) Here $\Hq(G)$ is the path $L_0L_1\cdots L_n$, so $c(G)=n+1$ and $\dim(\Hq(G))=\sdim(\Hq(G))=1$; moreover $\diam(\G(G))=n$, the distance from $e$ to a generator of $G$ (Lemma~\ref{lem:dist}). By Theorem~\ref{thm:dimcyc}(i), $\sdim(\G(G))=p^n-n=|G|-c(G)+1=|G|-n$. By Theorem~\ref{thm:dimcyc}(ii),(iii), $\dim(\G(G))=p^n-n=|G|-c(G)+1$ exactly in the listed cases, and $\dim(\G(G))=p^n-n-1=|G|-c(G)$ otherwise.
	\end{proof}
	
	\subsection{Groups attaining the lower bounds}
	Call a cyclic subgroup $E$ \emph{heavy} if $w(E)=|\gen(E)|\ge2$, i.e.\ $|E|\ge3$, and \emph{light} if $w(E)=|\gen(E)|=1$ (so the light cyclic subgroups are $1$ and the subgroups of order $2$). For cyclic subgroups $C,E$ put
	\[ \delta_C(E)=d_\Hq(C,E)\ \ (E\ne C),\qquad \delta_C(C)=2 . \]
	By Lemma~\ref{lem:dist}, $\delta_C(E)$ is the distance in $\G(G)$ from a fixed element of $\gen(C)$ to any \emph{other} element of $\gen(E)$.
	
	\begin{theorem}\label{thm:lowattain}
		Let $G$ be a finite group with $|G|\ge2$.
		\begin{enumerate}[label=(\roman*)]
			\item $\sdim(\G(G))=|G|-c(G)$ if and only if $G\cong\Z_p$ with $p$ an odd prime. For every other $G$, $\sdim(\G(G))\ge|G|-c(G)+1$.
			\item $\dim(\G(G))=|G|-c(G)$ if and only if for any two distinct cyclic subgroups $C,D$ of $G$ there is a heavy cyclic subgroup $E$ with $\delta_C(E)\ne\delta_D(E)$.
			\item If $\Hq(G)$ has two distinct vertices $C,D$ with $N_\Hq(C)=N_\Hq(D)$, then $\dim(\G(G))>|G|-c(G)$. In particular this holds if $G$ has two distinct maximal cyclic subgroups with the same set of maximal subgroups, for instance two distinct maximal cyclic subgroups of the same prime order.
		\end{enumerate}
	\end{theorem}
	\begin{proof}
		(i) If $G\cong\Z_p$ with $p$ odd, then $\G(G)=K_{1,p-1}$, $c(G)=2$, and $\sdim=p-2=|G|-c(G)$ by Theorem~\ref{thm:dimcyc}(i). Conversely, by Theorem~\ref{thm:sdim}, $\sdim(\G(G))=|G|-c(G)$ means $\alpha(Q(G))=c(G)$, i.e.\ $Q(G)$ has no edge. If $\diam(\Hq(G))\ge2$, let $C,D$ be at distance $\diam(\Hq(G))$ in $\Hq(G)$; then $d_\Hq(W,D)\le d_\Hq(C,D)$ for every $W$ and $C\not\sim D$, so (b1), (b2) of Lemma~\ref{lem:MMD} hold for $(C,D)$ and for $(D,C)$, and $CD$ is an edge of $Q(G)$. Hence $\diam(\Hq(G))\le1$, so $\Hq(G)$ is complete; being bipartite and connected (Lemmas~\ref{lem:bip},~\ref{lem:dist}) it is $K_2$, i.e.\ $c(G)=2$, which forces $G\cong\Z_p$. Finally $p\ne2$, since $\sdim(\G(\Z_2))=\sdim(K_2)=1>0$. In all other cases $\alpha(Q(G))\le c(G)-1$.
		
		(ii) By the proof of Theorem~\ref{thm:dim}, a resolving set $W$ with $|W|=|G|-c(G)$ omits exactly one element $x_C$ of each class $\gen(C)$, and any two such sets are mapped onto each other by twin automorphisms; so either all of them resolve or none does. Fix such a $W$. A pair containing an element of $W$ is resolved by that element. The class $\gen(E)$ meets $W$ if and only if $E$ is heavy, and then each $w\in W\cap\gen(E)$ satisfies $d(x_C,w)=\delta_C(E)$ for every $C$. Hence $W$ resolves if and only if for all $C\ne D$ some heavy $E$ has $\delta_C(E)\ne\delta_D(E)$.
		
		(iii) Since $\Hq(G)$ is bipartite, connected and $N(C)=N(D)\neq\emptyset$, we have $C\not\sim D$ and $d_\Hq(C,D)=2$. For $E\notin\{C,D\}$, a shortest $E$--$C$ path ends with an edge $FC$, where $F\in N(C)=N(D)$; so $d_\Hq(E,D)\le d_\Hq(E,C)$, and by symmetry equality holds. Moreover $\delta_C(C)=2=d_\Hq(D,C)=\delta_D(C)$, and likewise at $D$. So the condition in (ii) fails for $(C,D)$. A maximal cyclic subgroup has no cover, so its $\Hq$-neighbours are exactly its maximal subgroups, which gives the last assertion.
	\end{proof}
	
	\begin{corollary}\label{cor:lownot}
		$\dim(\G(G))>|G|-c(G)$ for the dihedral groups $D_n$ $(n\ge2)$, the generalized quaternion groups $\Q_{2^m}$, the elementary abelian groups of rank at least $2$, and for $S_3$, $S_4$, $A_4$, $A_5$.
	\end{corollary}
	\begin{proof}
		In $D_n$ the $n\ge2$ subgroups generated by reflections are maximal cyclic of order $2$. In $\Q_{2^m}=\langle x,y\rangle$ the $2^{m-2}\ge2$ subgroups $\langle x^iy\rangle$ are maximal cyclic of order $4$ and all have the same maximal subgroup $\langle x^{2^{m-2}}\rangle$. An elementary abelian $p$-group of rank $\ge2$ has at least $p+1\geq 3$ subgroups of order $p$, all maximal cyclic. In $S_3$ and $S_4$ the transpositions generate maximal cyclic subgroups of order $2$, and in $A_4$ and $A_5$ the subgroups of order $3$ (resp.\ $5$) are maximal cyclic. Apply Theorem~\ref{thm:lowattain}(iii).
	\end{proof}
	
	For cyclic groups the criterion of Theorem~\ref{thm:lowattain}(ii) can be decided completely. Recall that $c(\Z_n)=\tau(n)$, the number of divisors of $n$.
	
	\begin{theorem}\label{thm:lowcyc}
		Let $n\ge2$. Then $\dim(\G(\Z_n))=n-\tau(n)$ if and only if $n\ne 16$ and $n$ is not of any of the forms
		\[ 2,\qquad p^2,\qquad p^3,\qquad pq,\qquad p^2q,\qquad 8q, \]
		where $p,q$ denote distinct primes and, in the last form, $q$ is odd. In particular the lower bound is attained whenever $n$ has at least three distinct prime divisors, for $n=p$ and $n=p^a$ with $p$ odd and $a\ge4$, for $n=2^a$ with $a\ge5$, and for $n=p^aq^b$ unless $(a,b)\in\{(1,1),(1,2),(2,1)\}$ or $n=8q$.
	\end{theorem}
	\begin{proof}
		Write $n=p_1^{a_1}\cdots p_k^{a_k}$. By Lemma~\ref{lem:product} (or directly), identifying the subgroup of order $\prod p_i^{z_i}$ with $z=(z_1,\dots,z_k)$, $\Hq(\Z_n)$ is the grid $B=\prod_i[0,a_i]$ with $d_\Hq(z,z')=\|z-z'\|_1$. The light points are $0$ and, if some $p_i=2$, the unit vector $\varepsilon=e_i$; note $\|0-\varepsilon\|_1=1$. We check the criterion of Theorem~\ref{thm:lowattain}(ii) for distinct $x,y\in B$.
		
		\emph{Reduction to distance two.} If $x$ is heavy and $\|x-y\|_1\ne2$, then $E=x$ works, since $\delta_x(x)=2\ne\|x-y\|_1=\delta_y(x)$; similarly if $y$ is heavy. If both are light, then $\{x,y\}=\{0,\varepsilon\}$ and every heavy $z$ works, since $\|z-\varepsilon\|_1=\|z\|_1\pm1$; a heavy point exists unless $n=2$, and $n=2$ fails ($\G(\Z_2)=K_2$). So let $\|x-y\|_1=2$. Then $\delta_x(z)=\delta_y(z)=2$ for $z\in\{x,y\}$, and the pair is resolved if and only if the set $R=R(x,y)=\{z\in B:\|z-x\|_1\ne\|z-y\|_1\}$ contains a heavy point other than $x,y$, i.e.\ if and only if $R\setminus\{x,y\}\not\subseteq\{0,\varepsilon\}$. Since only the coordinates in which $x,y$ differ matter, there are three kinds of pairs (up to swapping $x,y$ and the indices $i\ne j$):
		\begin{itemize}
			\item[($\alpha$)] $y=x+2e_i$: \ $R=\{z: z_i\ne x_i+1\}$;
			\item[($\beta$)] $y=x+e_i+e_j$: \ $R=\{z: z_i\le x_i,\ z_j\le x_j\}\cup\{z: z_i> x_i,\ z_j> x_j\}$;
			\item[($\gamma$)] $x=w+e_i$, $y=w+e_j$: \ $R=\{z: z_i> w_i,\ z_j\le w_j\}\cup\{z: z_i\le w_i,\ z_j> w_j\}$.
		\end{itemize}
		(For instance in ($\beta$), $\|z-y\|_1-\|z-x\|_1$ is the sum of a term $+1$ or $-1$ according as $z_i\le x_i$ or $z_i>x_i$, and the same term for $j$.)
		
		\emph{A general observation.} In each case $R$ is a union of ``boxes'' $B'\times\prod_{l\in L}[0,a_l]$, where $L$ is the set of coordinates other than $i$ (resp.\ $i,j$). If $L\ne\emptyset$, then changing one coordinate $l\in L$ in both $x$ and $y$ to the same other value gives two further points $x',y'\in R$ with $\|x'-y'\|_1=\|x-y\|_1=2$; they cannot both be light, because the light points are at distance $1$. Hence the pair is resolved whenever $|R\setminus\{x,y\}|\ge3$, or $|R\setminus\{x,y\}|=2$ and $L\ne\emptyset$.
		
		\emph{$k\ge3$.} In ($\alpha$), $|R|=a_i\prod_{l\ne i}(a_l+1)\ge 2\cdot 4$. In ($\beta$), ($\gamma$), $|R|\ge2\prod_{l\ne i,j}(a_l+1)\ge4$, and $L\ne\emptyset$. Every pair is resolved.
		
		\emph{$k=1$.} Only ($\alpha$) occurs: $y=x+2$, and $R\setminus\{x,y\}=[0,a]\setminus\{x,x+1,x+2\}$. For $p$ odd this lies in $\{0\}$ if and only if $a=2$ ($x=0$) or $a=3$ ($x=1$); for $p=2$ it lies in $\{0,1\}$ if and only if $a\in\{2,3,4\}$ ($x=a-2$). This recovers Theorem~\ref{thm:dimcyc}: the forms $p^2$, $p^3$, $16$ (and $2$) fail, and $p$, $p^a$ ($p$ odd, $a\ge4$), $2^a$ ($a\ge5$) attain.
		
		\emph{$k=2$}, with exponents $a=a_1$, $b=a_2$. In ($\alpha$) with $i=1$, $|R|=a(b+1)\ge4$, with equality only when $a=2$, $b=1$; then $\{x_1,y_1\}=\{0,2\}$ and $R\setminus\{x,y\}$ contains a point with first coordinate $2$, which is heavy. So ($\alpha$) pairs are resolved.
		
		In ($\beta$) with $x=(s,t)$, $R$ is the union of the lower box $[0,s]\times[0,t]\ni 0,x$ and the upper box $[s+1,a]\times[t+1,b]\ni y$. The points of the upper box have both coordinates positive, so they are heavy; hence failure requires $(a,b)=(s+1,t+1)$ and $[0,s]\times[0,t]\setminus\{x\}\subseteq\{0,\varepsilon\}$, i.e.\ $(s+1)(t+1)\le3$. This gives exactly: $(s,t)=(0,0)$, $n=p_1p_2$; $(s,t)\in\{(1,0),(0,1)\}$, the extra point being $0$, $n=p^2q$; and $(s,t)=(2,0)$ (or $(0,2)$), the extra points being $0$ and $e_1$, which are both light if and only if $p_1=2$, i.e.\ $n=8q$.
		
		In ($\gamma$), $0\notin R$, so failure requires $R\setminus\{x,y\}\subseteq\{\varepsilon\}$. With $i=1$, $j=2$, the two boxes have sizes $(a-w_1)(w_2+1)$ and $(w_1+1)(b-w_2)$, each $\ge1$, with sum $\le3$. If both are $1$, then $a=b=1$ and $n=pq$. Otherwise, say the first box has $2$ points and $w_1=0$, $w_2=b-1$; either $a=2$, $b=1$, and the extra point $(2,0)$ is heavy, or $a=1$, $b=2$, and the extra point is $(1,0)=e_1$, light if and only if $p_1=2$, giving $n=2q^2$, which is of the form $p^2q$. So ($\gamma$) gives no new failures.
		
		Thus for $k=2$ the lower bound fails exactly when $n$ is of the form $pq$, $p^2q$ or $8q$.
	\end{proof}

	\section{Hamiltonicity}\label{sec:ham}
	
	A Hamiltonian cycle requires at least three vertices and minimum degree at least $2$, and a Hamiltonian path can have at most two vertices of degree $1$ (its endpoints). We first give a simple obstruction that applies to many non-nilpotent groups, including all dihedral groups, and then treat nilpotent groups.
	
	\subsection{Leaves and dihedral groups}
	
	\begin{lemma}\label{lem:leaf}
		Let $t\in G$ have prime order, and suppose that $\langle t\rangle$ is a maximal cyclic subgroup of $G$. Then $N(t)=\{e\}$. Consequently, if $|G|\ge3$ and $G$ has such an element, then $\G(G)$ is not Hamiltonian, and if $G$ has at least three such elements, then $\G(G)$ has no Hamiltonian path.
	\end{lemma}
	\begin{proof}
		By Lemma~\ref{lem:prime}, a neighbour $x\ne e$ of $t$ satisfies $t\in\langle x\rangle$ with $\ord x/\ord t$ prime, so $\langle t\rangle\subsetneq\langle x\rangle$, contradicting maximality; and $e\sim t$.
	\end{proof}
	
	\begin{corollary}\label{cor:dihham}
		For $n\ge2$, $\G(D_n)$ is not Hamiltonian, and for $n\ge3$ it has no Hamiltonian path; $\G(D_2)=K_{1,3}$ has no Hamiltonian path either.
	\end{corollary}
	\begin{proof}
		A cyclic subgroup containing a reflection is not generated by a rotation, so it is generated by a reflection and has order $2$. Hence the $n$ reflections generate maximal cyclic subgroups of order $2$, and Lemma~\ref{lem:leaf} applies.
	\end{proof}
	
	\subsection{Two tools}
	The blow-up structure reduces Hamiltonicity to a question about walks in $\Hq(G)$.
	
	\begin{lemma}\label{lem:walk}
		\begin{enumerate}[label=(\roman*)]
			\item $\G(G)$ has a Hamiltonian cycle if and only if $|G|\ge3$ and $\Hq(G)$ has a closed walk that visits every cyclic subgroup $C$ exactly $w(C)$ times.
			\item $\G(G)$ has a Hamiltonian path if and only if $\Hq(G)$ has a walk that visits every cyclic subgroup $C$ exactly $w(C)$ times.
		\end{enumerate}
	\end{lemma}
	\begin{proof}
		By Lemma~\ref{lem:blowup}, mapping each element to its class sends a Hamiltonian cycle (path) of $\G(G)$ to a closed walk (walk) of $\Hq(G)$ with the stated visit numbers. Conversely, given such a walk, assign to the $w(C)$ visits of $C$ the $w(C)$ distinct elements of $\gen(C)$; consecutive visits lie in adjacent classes, whose elements are pairwise adjacent by Lemma~\ref{lem:blowup}(ii).
	\end{proof}
	
	Let $X\subseteq G$ and $\Ev(X)=\{x:\Omega(\ord x)\text{ even}\}$ and $\Od(X)=\{x:\Omega(\ord x)\text{ odd}\}$. Define the invariant $s(X)=\sum_{x\in X}(-1)^{\Omega(\ord x)}=|\Ev(X)|-|\Od(X)|$ for $X\subseteq G$.
	
	\begin{lemma}\label{lem:obstr}
		\begin{enumerate}[label=(\roman*)]
			\item If $\G(G)$ is Hamiltonian then $s(G)=0$; if it has a Hamiltonian path then $|s(G)|\le1$.
			\item Let $\mathcal S$ be a set of cyclic subgroups of $G$, of total weight $w(\mathcal S)$, such that $\Hq(G)-\mathcal S$ has $m$ components. If $\G(G)$ is Hamiltonian then $m\le w(\mathcal S)$; if it has a Hamiltonian path then $m\le w(\mathcal S)+1$.
		\end{enumerate}
	\end{lemma}
	\begin{proof}
		(i) $\G(G)$ is bipartite with parts $\Ev(G),\Od(G)$ (Lemma~\ref{lem:bip}), and a Hamiltonian cycle (path) alternates between them. (ii) Deleting the $w(\mathcal S)$ elements of the classes in $\mathcal S$ leaves at least $m$ components, since there are no edges between classes in different components of $\Hq(G)-\mathcal S$; a Hamiltonian cycle (path) leaves at most $w(\mathcal S)$ (resp.\ $w(\mathcal S)+1$) components after deleting $w(\mathcal S)$ vertices.
	\end{proof}
	
	\subsection{Nilpotent groups}
	
	\begin{lemma}\label{lem:s}
		(a) If $G=P_1\times\cdots\times P_k$ is nilpotent, $s(G)=\prod_i s(P_i)$.
		(b) If $P\neq1$ is a $p$-group with $p$ odd, then $s(P)$ is odd and $s(P)\equiv2\pmod p$.
		(c) $s(\Z_2)=0$, and $s(P)\equiv2\pmod4$ for every $2$-group $P$ with $|P|\ge4$.
	\end{lemma}
	\begin{proof}
		Let $n_j$ be the number of elements of order $p^j$ in $P$. (a) $\ord{x_1\cdots x_k}=\prod\ord{x_i}$ for $x_i\in P_i$, so $\Omega$ is additive. (b) $s(P)\equiv|P|\pmod2$. By Fact~\ref{fact:group}(ii), $1+n_1\equiv0\pmod p$, and $p\mid n_j$ for $j\ge2$ because $n_j$ is a multiple of $\varphi(p^j)$; hence $s(P)\equiv1-n_1\equiv2\pmod p$. (c) By Fact~\ref{fact:group}(ii), $1+n_1$ is even and $1+n_1+n_2\equiv0\pmod4$; and $4\mid n_j$ for $j\ge3$. Hence $s(P)\equiv 1-n_1+n_2\equiv-2n_1\equiv2\pmod4$.
	\end{proof}
	
	For a $p$-group $P$ and a cyclic subgroup $C\le P$, let $c_P(C)$ be the number of cyclic subgroups of $P$ that contain $C$ as a maximal subgroup; by Lemma~\ref{lem:tree} these are the children of $C$ in the rooted tree $\Hq(P)$, and $N_p(P):=c_P(1)$ is the number of subgroups of order $p$. Let $G=K\times\langle z\rangle$ with $\ord z=2$ and $|K|$ odd. For a prime $p$ dividing $|K|$, let $P$ be the Sylow $p$-subgroup of $K$ and $m_p=|K|/|P|$, and consider the conditions
	\begin{itemize}
		\item[(H1)] $N_p(P)\le 2m_p$;
		\item[(H2)] $c_P(C)\le 2\varphi(|C|)\,m_p-1$ for every cyclic subgroup $1\ne C\le P$.
	\end{itemize}
	
	\begin{theorem}\label{thm:hamnec}
		Let $G$ be nilpotent with $|G|\ge3$. If $\G(G)$ is Hamiltonian, then $G\cong K\times\Z_2$ with $|K|>1$ odd, and (H1) and (H2) hold for every prime $p$ dividing $|K|$.
	\end{theorem}
	\begin{proof}
		By Lemma~\ref{lem:obstr}(i), $s(G)=0$, so by Lemma~\ref{lem:s} the Sylow $2$-subgroup is $\Z_2$; and $K\ne1$ since $|G|\ge3$. Fix $p$, and write $G=P\times R$ with $R$ the product of the other Sylow subgroups of $G$, so $|R|=2m_p$. By Lemma~\ref{lem:product}, $\Hq(G)=\Hq(P)\,\square\,X$ with $X=\Hq(R)$ connected and of total weight $|R|$, and $w(C,x)=\varphi(|C|)\,w_X(x)$. For a cyclic $C\le P$ let $\mathcal S_C=\{C\}\times X$, of weight $2\varphi(|C|)m_p$. Since $\Hq(P)$ is a tree, $\Hq(P)-C$ has $c_P(C)$ components if $C=1$ and $c_P(C)+1$ components otherwise, and $\Hq(G)-\mathcal S_C$ is the product of these with $X$. Lemma~\ref{lem:obstr}(ii) gives (H1) for $C=1$ and (H2) for $C\ne1$.
	\end{proof}
	
	The conditions are also sufficient when all Sylow subgroups of $K$ except one are cyclic.
	
	\begin{theorem}\label{thm:hamsuf}
		Let $K=P\times R$, where $P\ne1$ is a $p$-group ($p$ odd) and $R$ is a cyclic group of odd order prime to $p$ ($R=1$ allowed), and let $G=K\times\Z_2$. Then $\G(G)$ is Hamiltonian if and only if
		\[ N_p(P)\le2|R|\qquad\text{and}\qquad c_P(C)\le 2\varphi(|C|)\,|R|-1\ \text{ for every cyclic }1\neq C\le P. \]
	\end{theorem}
	\begin{proof}
		Necessity is Theorem~\ref{thm:hamnec} for the prime $p$ (here $m_p=|R|$). For sufficiency we induct on the number of prime divisors of $|R|$. Let $X=\Hq(R\times\Z_2)$, with weights $w_X$ of total $\ell=2|R|$. We first need a closed walk $W$ in $X$ of length $\ell$ visiting each $x$ exactly $w_X(x)$ times. If $R=1$, then $X$ is the edge $\{1,\Z_2\}$ and $W=(1,\Z_2)$. If $R\ne1$, write $R=Q\times R'$ with $Q\ne1$ a cyclic $q$-group; since $N_q(Q)=1$ and $c_Q(C)\le1$, the hypotheses hold for $(Q,R')$, so $\G(R\times\Z_2)$ is Hamiltonian by induction and $W$ exists by Lemma~\ref{lem:walk}(i).
		
		By Lemma~\ref{lem:product}, $\Hq(G)=\Hq(P)\,\square\,X$ with $w(C,x)=\varphi(|C|)w_X(x)$. For a vertex $C$ of the tree $\Hq(P)$ let $T(C)$ be the set of $C$ and its descendants, and let $W^m$ denote the closed walk that runs $m$ times around $W$. A \emph{step} of $W$ is a pair $x\to x'$ of consecutive vertices.
		
		\emph{Claim.} For every $C\neq1$ in $\Hq(P)$ and every step $x\to x'$ of $W$, the subgraph of $\Hq(G)$ induced by $T(C)\times X$ has a closed walk $\Omega_C$ visiting every vertex exactly its weight times and containing the step $(C,x)\to(C,x')$.
		
		We prove the claim by induction on the height of $C$ in the tree, and simultaneously build the required walk on $\Hq(G)=T(1)\times X$. Start with the closed walk $\{C\}\times W^{\varphi(|C|)}$; it visits $(C,y)$ exactly $\varphi(|C|)w_X(y)=w(C,y)$ times and has $\varphi(|C|)\,\ell$ steps. If $C\ne1$, reserve one occurrence of the step $(C,x)\to(C,x')$ (it occurs $\varphi(|C|)\ge2$ times). The number of remaining steps is $\ell$ if $C=1$, and $\varphi(|C|)\ell-1$ if $C\ne1$; by the hypotheses this is at least the number $c_P(C)$ of children of $C$. Assign to the children $D$ of $C$ distinct unreserved steps $(C,y_D)\to(C,y'_D)$. By induction, $\Omega_D$ exists with the step $(D,y_D)\to(D,y'_D)$; deleting this step from $\Omega_D$ and reversing gives a walk $\Pi_D$ from $(D,y_D)$ to $(D,y'_D)$ visiting every vertex of $T(D)\times X$ exactly its weight times. Replace the step $(C,y_D)\to(C,y'_D)$ by $(C,y_D)\,\Pi_D\,(C,y'_D)$; this is a walk because $(C,y)\sim(D,y)$ in $\Hq(P)\,\square\,X$. Doing this for all children gives a closed walk $\Omega_C$ with the required visit numbers, which still contains the reserved step. For $C=1$ we obtain a closed walk of $\Hq(G)$ visiting each vertex exactly its weight times, of length $|G|\ge6$, and Lemma~\ref{lem:walk}(i) finishes the proof.
	\end{proof}
	
	\begin{corollary}\label{cor:hamab}
		\begin{enumerate}[label=(\roman*)]
			\item $\G(\Z_n)$ is Hamiltonian if and only if $n=2m$ with $m>1$ odd.
			\item If $K$ is a $p$-group ($p$ odd), then $\G(K\times\Z_2)$ is Hamiltonian if and only if $K$ is cyclic.
			\item Let $K=P\times R$, where $P$ is an abelian $p$-group of rank $r$ ($p$ odd) and $R$ is cyclic of odd order prime to $p$. Then $\G(K\times\Z_2)$ is Hamiltonian if and only if
			\[ \frac{p^r-1}{p-1}\ \le\ 2|R| . \]
		\end{enumerate}
	\end{corollary}
	\begin{proof}
		(i) Necessity is Theorem~\ref{thm:hamnec}; if $m>1$ is odd, apply Theorem~\ref{thm:hamsuf} with $P$ a nontrivial cyclic Sylow subgroup of $\Z_m$, for which $N_p(P)=1$ and $c_P(C)\le1$.
		(ii) Take $R=1$: the first condition reads $N_p(K)\le2$. By Fact~\ref{fact:group}(iv), $N_p(K)\equiv1\pmod p$ with $p\ge3$, so $N_p(K)=1$ and $K$ is cyclic by Fact~\ref{fact:group}(iii); conversely, cyclic $K$ is covered by (i).
		(iii) Here $N_p(P)=(p^r-1)/(p-1)$, and we show that the second condition of Theorem~\ref{thm:hamsuf} follows from the first. Let $C=\langle g\rangle\ne1$. Each child $D$ of $C$ has a generator $y$ with $y^p=g$ (if $y_0$ generates $D$, then $y_0^p=g^j$ with $p\nmid j$, and $y=y_0^{j'}$ with $jj'\equiv1\pmod{|C|}$ works), and $D$ contains exactly $p$ such elements, namely $y^t$ with $t\equiv1\pmod{|C|}$. In the abelian group $P$, the solutions of $y^p=g$ form a coset of $\{u:u^p=1\}$, which has order $p^r$. Hence $c_P(C)\le p^{r-1}\le N_p(P)\le2|R|<2\varphi(|C|)|R|-1$, since $\varphi(|C|)\ge2$.
	\end{proof}

	\subsection{Hamiltonian paths}
	
	\begin{theorem}\label{thm:hampath}
		Let $G$ be nilpotent with $|G|\ge2$.
		\begin{enumerate}[label=(\roman*)]
			\item If $\G(G)$ has a Hamiltonian path, then $G\cong\Z_3$ or $G\cong K\times\Z_2$ with $|K|$ odd.
			\item If $G=K\times\Z_2$ and $\G(G)$ has a Hamiltonian path, then (H1) holds for every prime $p$ dividing $|K|$.
			\item In the situation of Corollary~\ref{cor:hamab}(iii), $\G(K\times\Z_2)$ has a Hamiltonian path if and only if it is Hamiltonian. In particular, $\G(\Z_n)$ has a Hamiltonian path if and only if $n\in\{2,3\}$ or $n=2m$ with $m>1$ odd.
		\end{enumerate}
	\end{theorem}
	\begin{proof}
		(i) By Lemma~\ref{lem:obstr}(i), $|s(G)|\le1$. Write $G=P_2\times K$ with $P_2$ the Sylow $2$-subgroup. By Lemma~\ref{lem:s}, $s(K)$ is odd, and $|s(P_2)|\ge2$ if $|P_2|\ge4$; hence $|P_2|\le2$. Suppose $G=K$ has odd order. Then $|s(P)|=1$ for every Sylow subgroup $P\neq1$; as $s(P)\equiv2\pmod p$, this forces $s(P)=-1$ and $p=3$. So $G$ is a nontrivial $3$-group, and $\Hq(G)$ is a tree in which the root $1$ has weight $1$ (Lemma~\ref{lem:tree}). By Lemma~\ref{lem:walk}(ii) there is a walk visiting each vertex its weight times; it visits $1$ once, so it meets at most two components of $\Hq(G)-1$, whence $N_3(G)\le2$ and $G\cong\Z_{3^n}$ by Fact~\ref{fact:group}(iii),(iv). If $n\ge2$, the vertex $L_n$ of order $3^n$ has $L_{n-1}$ as its only neighbour, so consecutive visits to $L_n$ are separated by visits to $L_{n-1}$ and $\varphi(3^n)\le\varphi(3^{n-1})+1$, which is false. So $G\cong\Z_3$, and indeed $\G(\Z_3)=K_{1,2}$ is a path.
		
		(ii) As in the proof of Theorem~\ref{thm:hamnec}, Lemma~\ref{lem:obstr}(ii) gives $N_p(P)\le2m_p+1$. Equality is impossible, since $N_p(P)\equiv1\pmod p$ by Fact~\ref{fact:group}(iv), while $p\nmid 2m_p$.
		
		(iii) By (ii) and Corollary~\ref{cor:hamab}(iii), a Hamiltonian path forces a Hamiltonian cycle; the converse is clear. For $\Z_n$, combine (i), (iii) and Corollary~\ref{cor:hamab}(i): the case $G=\Z_2$ gives the path $K_2$, and $G=\Z_3$ gives $K_{1,2}$.
	\end{proof}
	
	\begin{problem}\label{prob:ham}
		Let $K$ be a nilpotent group of odd order. Are conditions (H1) and (H2), for all primes $p$ dividing $|K|$, sufficient for $\G(K\times\Z_2)$ to be Hamiltonian? Does (H1) imply (H2), as it does for abelian Sylow subgroups (proof of Corollary~\ref{cor:hamab}(iii))? Does $\G(K\times\Z_2)$ always have a Hamiltonian cycle when it has a Hamiltonian path?
	\end{problem}
	\section{Genus}\label{sec:genus}
	
	We recall the topological notions used here (see, e.g.,~\cite{Ringel}). A \emph{closed surface} is a compact connected $2$-manifold without boundary; $S_g$ denotes the orientable closed surface of genus $g$, i.e.\ the sphere with $g$ handles. An \emph{embedding} of a graph $X$ in a surface $S$ is a drawing of $X$ on $S$ in which edges meet only at common end vertices. The connected components of $S\setminus X$ are the \emph{faces} of the embedding. The embedding is \emph{cellular} (or a \emph{$2$-cell embedding}) if every face is homeomorphic to an open disc. In that case each face $f$ is bounded by a closed walk of $X$, its \emph{boundary walk}, whose length $\ell_f$ is the \emph{length} of $f$. Every edge occurs exactly twice in the boundary walks, and on an orientable surface it is traversed once in each direction when all faces are oriented consistently; hence $\sum_f\ell_f=2E$. Cellular embeddings of a connected graph in orientable surfaces correspond to \emph{rotation systems}, that is, a cyclic order of the neighbours at each vertex. The faces are traced by the rule that, after traversing the edge $uv$, one leaves $v$ along the successor of $u$ in the cyclic order at $v$. The \emph{genus} $\genus(X)$ is the least $g$ such that $X$ embeds in $S_g$; $X$ is \emph{planar} if $\genus(X)=0$ and \emph{toroidal} if $\genus(X)=1$. By Fact~\ref{fact:graph}(ii), a minimum-genus embedding is cellular and satisfies Euler's formula $V-E+F=2-2g$.
	
	Throughout this section
	\[ \Pi=\{1,2,3,4,6,8,9\}\cup\{p\ge5:\ p\text{ prime}\} \]
	is the set of \emph{admissible orders}. We proceed in three steps of increasing difficulty. Planarity needs only Kuratowski's theorem (Section~\ref{ss:planar}); toroidality needs, in addition, Euler's formula and a short counting argument (Sections~\ref{ss:comp}--\ref{ss:torus}); the genus-two case needs, in addition, a few lower bounds of the form $\genus\ge3$, which we derive from the same counting argument, Euler's formula and some elementary group theory (Section~\ref{ss:genus2}). Since $\G(H)$ is an induced subgraph of $\G(G)$ for $H\le G$ (Lemma~\ref{lem:blowup}), the genus of $\G(G)$ is at least the genus of $\G(\langle x\rangle)$ for every $x\in G$.
	
	\subsection{Element orders}
	
	\begin{lemma}\label{lem:orders}
		Let $G$ have an element of order $n\notin\Pi$. Then $\G(G)$ contains $K_{3,3}$, and $\genus(\G(G))\ge2$. Moreover $\genus(\G(G))\ge3$.
	\end{lemma}
	\begin{proof}
		If $n\notin\Pi$, then $n$ is divisible by one of the following numbers $d$, and $\G(\Z_d)\le\G(G)$:
		\begin{itemize}
			\item $d=pq\ne6$ ($p<q$ primes): $\G(\Z_{pq})=K_{p+q-2,\,(p-1)(q-1)+1}\supseteq K_{5,5}$, genus $\ge3$ by Fact~\ref{fact:graph}(iv);
			\item $d=p^2$, $p\ge5$: the elements of orders $p$ and $p^2$ form $K_{p-1,p^2-p}\supseteq K_{4,20}$, genus $\ge9$;
			\item $d=16$: the elements of orders $8$ and $16$ form $K_{4,8}$, genus $3$;
			\item $d=27$: the elements of orders $9$ and $27$ form $K_{6,18}$, genus $16$;
			\item $d=18$: $V=18$, $E=69$, genus $\ge\lceil 37/4\rceil=10$ by Fact~\ref{fact:graph}(v);
			\item $d=12$: $V=12$, $E=27$, genus $\ge\lceil 7/4\rceil=2$ by Fact~\ref{fact:graph}(v), and $\G(\Z_{12})\supseteq\G(\Z_6)=K_{3,3}$.
		\end{itemize}
		Indeed, if two primes divide $n$ and they are not $\{2,3\}$, then some $pq\ne6$ divides $n$; if $n=2^a3^b$ with $a,b\ge1$ and $n\ne6$, then $12\mid n$ or $18\mid n$; and a prime power outside $\Pi$ is divisible by $p^2$ ($p\ge5$), $16$ or $27$. This proves the first two assertions. For the last one only $d=12$ remains, and $\genus(\G(\Z_{12}))\ge3$ is Lemma~\ref{lem:Z12}.
	\end{proof}
	
	\subsection{Planar graphs}\label{ss:planar}
	
	\begin{theorem}[Planarity]\label{thm:planar}
		$\G(G)$ is planar if and only if every element order of $G$ lies in $\{1,2,3,4,8,9\}\cup\{p\ge5\text{ prime}\}$, i.e.\ all element orders lie in $\Pi$ and $G$ has no element of order $6$.
	\end{theorem}
	\begin{proof}
		If $G$ has an element of order $6$, then $\G(G)\supseteq\G(\Z_6)=K_{3,3}$; if it has an element of order outside $\Pi$, then $\G(G)\supseteq K_{3,3}$ by Lemma~\ref{lem:orders}. In both cases $\G(G)$ is not planar by Kuratowski's theorem.
		
		Conversely, suppose that every element order lies in $\{1,2,3,4,8,9\}\cup\{p\ge5\}$. Since $x\sim y$ requires $\ord y=q\cdot\ord x$ with $q$ prime, every edge of $\G(G)$ is of one of the forms: $e$--(element of prime order); involution--(element of order $4$ above it); (order $4$)--(order $8$ above it); (order $3$)--(order $9$ above it). Hence $\G(G)$ is obtained from the star on $e$ and the elements of prime order by attaching, at single vertices,
		\begin{itemize}
			\item for each cyclic subgroup $C$ of order $4$, with involution $t$, the graph $K_{2,1+4r}$ with parts $\gen(C)$ and $\{t\}\cup\{\text{the $4r$ elements of order $8$ above }C\}$, attached at $t$;
			\item for each subgroup $U$ of order $3$, the graph $K_{2,1+6r}$ with parts $\gen(U)$ and $\{e\}\cup\{\text{the $6r$ elements of order $9$ above }U\}$, attached at $e$.
		\end{itemize}
		All blocks are edges or graphs $K_{2,s}$, which are planar, so $\G(G)$ is planar by Fact~\ref{fact:graph}(iii).
	\end{proof}
	
	\subsection{Decomposition into components}\label{ss:comp}
	From now on assume that all element orders of $G$ lie in $\Pi$. Let $\Bc(G)$ be the bipartite graph whose vertices are the involutions and the subgroups of order $3$, with $t\sim U$ if and only if $t$ centralizes $U$; its edges correspond bijectively to the cyclic subgroups $\langle t\rangle\times U$ of order $6$. For a component $K$ of $\Bc(G)$ with an edge, let $V_K$ consist of $e$, the involutions of $K$, the elements of the subgroups in $K$, the elements of order $6$ generating edges of $K$, and the elements of order $9$ whose cube lies in a subgroup of $K$; put $\G_K=\G(G)[V_K]$.
	
	\begin{lemma}\label{lem:decomp}
		$\genus(\G(G))=\sum_K\genus(\G_K)$, the sum over the components $K$ of $\Bc(G)$ with an edge.
	\end{lemma}
	\begin{proof}
		Now the edges of $\G(G)$ are those listed in the proof of Theorem~\ref{thm:planar}, together with (order $3$)--(order $6$ over its subgroup) and involution--(order $6$ element whose cube it is). Hence distinct $\G_K$ meet only in $e$, and the rest of $\G(G)$ consists of the planar pieces of the proof of Theorem~\ref{thm:planar} (leaves of prime order $\ge5$, the graphs $K_{2,1+4r}$, and the graphs $K_{2,1+6r}$ of the order-$3$ subgroups outside all components), each attached at a single cut vertex. Fact~\ref{fact:graph}(iii) gives the claim.
	\end{proof}
	
	\begin{lemma}\label{lem:Bdeg}
		In $\Bc(G)$ every subgroup of order $3$ of positive degree has odd degree, and every involution of positive degree has degree $\equiv1\pmod3$.
	\end{lemma}
	\begin{proof}
		The degree of $U$ is the number of involutions in $C_G(U)$; if it is positive, $C_G(U)$ has even order and the number is odd by Fact~\ref{fact:group}(ii). The degree of $t$ is half the number of elements of order $3$ in $C_G(t)$; if it is positive, $3\mid|C_G(t)|$, that number is $\equiv2\pmod3$ by Fact~\ref{fact:group}(ii), and its half is $\equiv1\pmod3$.
	\end{proof}
	
	\begin{lemma}\label{lem:no5}
		No subgroup of order $3$ has degree $5$ in $\Bc(G)$.
	\end{lemma}
	\begin{proof}
		Let $U=\langle u\rangle$ have order $3$ and $H=C_G(U)$, so $u$ is central in $H$. If $h\in H$ has order $n$ prime to $3$, then $hu$ has order $3n\in\Pi$, which forces $n\in\{1,2\}$. Hence $H$ is a $\{2,3\}$-group (Fact~\ref{fact:group}(i)) without elements of order $4$, so its Sylow $2$-subgroups are elementary abelian, of order $2^r$ say, and every involution of $H$ lies in one of them. Suppose $H$ has exactly five involutions. Then $2^r-1\le5$, so $r\le2$, and the number $n_2$ of Sylow $2$-subgroups of $H$ is a power of $3$, as it divides $|H|/2^r$. If $r=1$, the involutions correspond bijectively to the Sylow $2$-subgroups, so $5=n_2$ is a power of $3$, which is absurd. If $r=2$ and $n_2=1$, there are only three involutions. If $r=2$ and $n_2\ge3$, three distinct Sylow subgroups $\cong\Z_2^2$ meet pairwise in at most one involution, so together they contain at least $3\cdot3-3=6$ involutions, again a contradiction.
	\end{proof}
	
	For $m\ge1$ let $S_U(m)$ be the graph on $e,w,W,t_i,x_i,x_i'$ ($1\le i\le m$) with edges $e w, eW, et_i, t_ix_i,t_ix_i'$ and $w x_i,wx_i',Wx_i,Wx_i'$ (this is $\G_K$ when $K$ is a star centred at $U=\{1,w,W\}$ with $m$ involutions $t_i$, $x_i,x_i'$ generating $\langle t_i\rangle\times U$, and no elements of order $9$). For $k\ge1$ let $S_t(k)$ be the graph on $e,t,w_i,W_i,x_i,x_i'$ with edges $et$, $ew_i,eW_i$, $tx_i,tx_i'$, $w_ix_i,w_ix_i',W_ix_i,W_ix_i'$ (the star centred at an involution $t$ with $k$ subgroups $U_i=\{1,w_i,W_i\}$). Note $S_U(1)=S_t(1)=K_{3,3}=\G(\Z_6)$.
	
	\begin{lemma}[Shape of the components]\label{lem:shape}
		Let $K$ be a component of $\Bc(G)$ with at least two edges. Then exactly one of the following holds:
		\begin{enumerate}[label=(\alph*)]
			\item $K$ is a star centred at a subgroup $U$ of order $3$ with $m$ involutions, where $m\ge3$ is odd; then $\G_K$ consists of $S_U(m)$ together with the elements of order $9$ above $U$, each adjacent exactly to the two generators of $U$;
			\item $K$ is a star centred at an involution $t$ with $k$ subgroups of order $3$, where $k\ge4$ and $k\equiv1\pmod3$; then $\G_K$ consists of $S_t(k)$ together with the elements of order $9$ above the $U_i$, each adjacent exactly to the two generators of its $U_i$;
			\item $K$ is not a star; then $K$ contains an involution $t$ of degree at least $4$ adjacent to a subgroup of order $3$ of degree at least $3$, and in particular $\G_K\supseteq S_t(4)$.
		\end{enumerate}
	\end{lemma}
	\begin{proof}
		In (a) and (b) the degree restrictions are Lemma~\ref{lem:Bdeg}, and the description of $\G_K$ follows from the list of edges in the proof of Lemma~\ref{lem:decomp}: an element of order $9$ has the two generators of its cube's subgroup as its only neighbours. If $K$ is not a star, then, being connected and bipartite, it contains a path $a\,b\,c\,d$; the adjacent vertices $b,c$ both have degree at least $2$, so by Lemma~\ref{lem:Bdeg} one of them is an involution $t$ of degree at least $4$ and the other a subgroup $U_1$ of order $3$ of degree at least $3$. The vertices $e$, $t$, the generators of four subgroups of order $3$ centralized by $t$ and the generators of the corresponding cyclic subgroups of order $6$ span a copy of $S_t(4)$.
	\end{proof}
	
	\subsection{Lower bounds by counting}
	
	Let $X$ be a connected bipartite graph with parts $A,B$ and minimum degree at least $3$, cellularly embedded in the orientable surface $S_g$. Inside each face of length $\ell\ge6$ insert a \emph{dummy} vertex joined to the $\ell/2$ corners of the face lying in $B$; every face is now a quadrilateral $b\,a\,b'\,a'$. Draw in each such face the diagonal $bb'$, and let $M$ be the resulting graph on $B$.
	
	\begin{lemma}\label{lem:balance}
		\begin{enumerate}[label=(\roman*)]
			\item $M$ is cellularly embedded in $S_g$; its faces correspond to $A\cup\{\text{dummies}\}$, the face of $a\in A$ being a polygon listing $N(a)$ in rotation order.
			\item If the faces of $M$ are oriented by the orientation of $S_g$, then for all distinct $b,b'\in B$, the number of face sides traversed $b\to b'$ equals the number traversed $b'\to b$.
			\item The total number of sides of dummy faces is at most $\tfrac32\,(4V-2E-8+8g)$, where $V,E$ refer to $X$.
		\end{enumerate}
	\end{lemma}
	\begin{proof}
		(i) Deleting the original edges, the union of the quadrilaterals at a vertex $a$ is a disk bounded by the diagonals of these quadrilaterals, in the rotation order at $a$. (ii) Each edge of $M$ borders two face-sides which an orientation traverses in opposite directions. (iii) A face of length $\ell\ge6$ contributes $\ell/2\le\tfrac32(\ell-4)$ dummy sides, and $\sum_f(\ell_f-4)=2E-4F=4V-2E-8+8g$ by Fact~\ref{fact:graph}(ii).
	\end{proof}
	
	\begin{lemma}\label{lem:starlow}
		\begin{enumerate}[label=(\roman*)]
			\item $\genus(S_U(m))\ge\lceil(4m+1)/12\rceil$; in particular $\genus(S_U(3))\ge2$ and $\genus(S_U(m))\ge3$ for $m\ge7$.
			\item $\genus(S_t(k))\ge\lceil k/3\rceil$; in particular $\genus(S_t(4))\ge2$ and $\genus(S_t(k))\ge3$ for $k\ge7$.
		\end{enumerate}
	\end{lemma}
	\begin{proof}
		Both $S_U(m)$ and $S_t(k)$ are bipartite with minimum degree $3$. Fix a cellular embedding of minimum genus $g$ (Fact~\ref{fact:graph}(ii)) and apply Lemma~\ref{lem:balance} with $B$ the side of odd $\Omega$ (the elements of orders $2$ and $3$). In both cases $N(e)=B$, so the face of $M$ at $e$ is a polygon through all vertices of $B$, each visited once; we call it the polygon $e$.
		
		(i) For $S_U(m)$: $B=\{w,W,t_1,\dots,t_m\}$, $V=3+3m$, $E=2+7m$, and the faces of $M$ are the polygon $e$, triangles $x_i,x_i'$ on $\{w,W,t_i\}$, and dummies. A pair containing $t_i$ receives an even total from the triangles (only $x_i,x_i'$ contain $t_i$), so each of the at least $m+1$ sides of the polygon $e$ that contain some $t_i$ needs a dummy side on the same pair. Lemma~\ref{lem:balance}(iii) gives $12g-3m\ge m+1$.
		
		(ii) For $S_t(k)$: $B=\{t,w_i,W_i\}$, $V=4k+2$, $E=8k+1$, and the faces of $M$ are the polygon $e$ (with $2k+1$ sides), the $2k$ triangles $x_i,x_i'$ on $\{t,w_i,W_i\}$ (with $6k$ sides), and the dummy faces, with $D$ sides in total. Every edge of $M$ is a side of exactly two faces, so the $8k+1$ sides of non-dummy faces are partitioned into pairs of sides of a common edge; let $P$ be the number of such pairs consisting of two non-dummy sides. Each remaining non-dummy side is paired with a dummy side, so $D\ge 8k+1-2P$. We bound $P$. A pair of sides of two triangles lies on a pair of vertices common to both, so only $x_i$ and $x_i'$ can be paired with each other. They share at most one side: if they shared two sides, these would meet at a vertex $v\in B$, and the two triangles would be the only faces of $M$ around $v$, whereas the polygon $e$ also passes through $v$. The polygon $e$ is not paired with itself (its vertices are distinct), and a side of $e$ can be paired with a triangle only if both its ends lie in a common triangle, i.e.\ if it contains $t$ (two sides) or is of the form $w_iW_i$ (at most $k$ sides). Hence $P\le k+(k+2)$ and $D\ge 4k-3$. On the other hand Lemma~\ref{lem:balance}(iii) gives $D\le\frac32(4V-2E-8+8g)=12g-3$. So $12g\ge4k$, i.e.\ $g\ge k/3$.
	\end{proof}
	
	\subsection{Embeddings}
	
	\begin{lemma}\label{lem:single}
		If $K$ is a single edge then $\genus(\G_K)=1$.
	\end{lemma}
	\begin{proof}
		$\G_K$ is $K_{3,3}$ together with vertices of order $9$, each adjacent exactly to the two generators $w,W$ of $U$. Embed $K_{3,3}$ in the torus with three hexagonal faces; each hexagon is a Hamiltonian cycle and so contains $w$ and $W$, and the order-$9$ vertices are drawn as chords $w$--$y$--$W$ in one hexagon. $K_{3,3}$ is non-planar.
	\end{proof}
	
	\begin{lemma}\label{lem:starup}
		\begin{enumerate}[label=(\roman*)]
			\item $\genus(S_U(3))=2$, and this remains true after adding any number of vertices adjacent exactly to $w$ and $W$.
			\item $\genus(S_t(4))=2$, and this remains true after adding, for each $i$, any number of vertices adjacent exactly to $w_i$ and $W_i$.
		\end{enumerate}
	\end{lemma}
	\begin{proof}
		Table~\ref{tab:emb2} gives rotation systems of $S_U(3)$ ($V=12$, $E=23$, $F=9$) and of $S_t(4)$ ($V=18$, $E=33$, $F=13$); in each case every edge occurs in exactly two of the listed face walks, once in each direction, and $V-E+F=-2$, so these are embeddings in the orientable surface of genus $2$. With Lemma~\ref{lem:starlow}, $\genus(S_U(3))=\genus(S_t(4))=2$. In the embedding of $S_U(3)$ the quadrilateral $w\,x_1\,W\,x_2$ is a face, and in the embedding of $S_t(4)$ each quadrilateral $w_i\,x_i\,W_i\,x_i'$ is a face. Any number of new vertices adjacent exactly to the two generators of the corresponding subgroup can be drawn inside such a face as pairwise disjoint paths of length $2$ between its two opposite corners; this does not change the genus.
	\end{proof}
	
	\begin{table}[ht]
		\centering\small
		\caption{Embeddings of $S_U(3)$ and $S_t(4)$ in the orientable surface of genus $2$ (rotation at each vertex, and the resulting face boundary walks).}\label{tab:emb2}
		\begin{tabular}{@{}p{0.47\textwidth}p{0.47\textwidth}@{}}
			\toprule
			$S_U(3)$: rotations & faces\\\midrule
			\raggedright $e:(w,t_3,t_2,t_1,W)$,\ $w:(e,x_2,x_1,x_2',x_3',x_1',x_3)$,\ $W:(x_3,x_1,x_2,x_3',x_2',e,x_1')$,\ $t_1:(x_1,x_1',e)$,\ $t_2:(e,x_2,x_2')$,\ $t_3:(x_3',e,x_3)$,\ $x_1:(t_1,w,W)$,\ $x_1':(w,W,t_1)$,\ $x_2:(W,w,t_2)$,\ $x_2':(W,w,t_2)$,\ $x_3:(w,W,t_3)$,\ $x_3':(t_3,w,W)$ & \raggedright $e\,w\,x_2\,t_2\,x_2'\,W$;\ \ $e\,t_3\,x_3\,w$;\ \ $e\,t_2\,x_2\,W\,x_3'\,t_3$;\ \ $e\,t_1\,x_1\,w\,x_2'\,t_2$;\ \ $e\,W\,x_1'\,t_1$;\ \ $w\,x_1\,W\,x_2$;\ \ $w\,x_3'\,W\,x_2'$;\ \ $w\,x_1'\,W\,x_3\,t_3\,x_3'$;\ \ $w\,x_3\,W\,x_1\,t_1\,x_1'$\tabularnewline\midrule
			$S_t(4)$: rotations & faces\\\midrule
			\raggedright $e:(W_1,W_4,W_3,w_2,w_1,W_2,t,w_4,w_3)$,\ $t:(e,x_2',x_3',x_4,x_1',x_2,x_1,x_3,x_4')$,\ $W_1:(x_1,x_1',e)$,\ $W_2:(x_2',e,x_2)$,\ $W_3:(x_3,x_3',e)$,\ $W_4:(x_4,x_4',e)$,\ $w_1:(e,x_1',x_1)$,\ $w_2:(e,x_2',x_2)$,\ $w_3:(e,x_3',x_3)$,\ $w_4:(x_4,e,x_4')$,\ $x_1:(W_1,t,w_1)$,\ $x_1':(t,W_1,w_1)$,\ $x_2:(w_2,W_2,t)$,\ $x_2':(t,W_2,w_2)$,\ $x_3:(t,w_3,W_3)$,\ $x_3':(w_3,t,W_3)$,\ $x_4:(w_4,W_4,t)$,\ $x_4':(w_4,t,W_4)$ & \raggedright $e\,W_1\,x_1\,t\,x_3\,w_3$;\ \ $e\,W_4\,x_4\,t\,x_1'\,W_1$;\ \ $e\,W_3\,x_3\,t\,x_4'\,W_4$;\ \ $e\,w_2\,x_2'\,t\,x_3'\,W_3$;\ \ $e\,w_1\,x_1'\,t\,x_2\,w_2$;\ \ $e\,W_2\,x_2\,t\,x_1\,w_1$;\ \ $e\,t\,x_2'\,W_2$;\ \ $e\,w_4\,x_4'\,t$;\ \ $e\,w_3\,x_3'\,t\,x_4\,w_4$;\ \ $w_1\,x_1\,W_1\,x_1'$;\ \ $w_2\,x_2\,W_2\,x_2'$;\ \ $w_3\,x_3\,W_3\,x_3'$;\ \ $w_4\,x_4\,W_4\,x_4'$\tabularnewline
			\bottomrule
		\end{tabular}
	\end{table}
	
	\subsection{Toroidal graphs}\label{ss:torus}
	
	\begin{theorem}[Toroidality]\label{thm:torus}
		$\G(G)$ is toroidal if and only if every element order of $G$ lies in $\Pi$ and $G$ has exactly one cyclic subgroup of order $6$.
	\end{theorem}
	\begin{proof}
		If some element order lies outside $\Pi$, then $\genus(\G(G))\ge2$ by Lemma~\ref{lem:orders}. Otherwise $\genus(\G(G))=\sum_K\genus(\G_K)$ by Lemma~\ref{lem:decomp}, the sum over the components $K$ of $\Bc(G)$ with an edge, whose edges correspond to the cyclic subgroups of order $6$. Each $\G_K$ contains $K_{3,3}=\G(\Z_6)$, so has genus at least $1$. If $K$ has at least two edges, Lemma~\ref{lem:shape} shows that $\G_K$ contains $S_U(3)$ or $S_t(4)$ , so $\genus(\G_K)\ge2$ by Lemma~\ref{lem:starlow}. If $K$ is a single edge, $\genus(\G_K)=1$ by Lemma~\ref{lem:single}. Hence, together with Theorem~\ref{thm:planar}, $\genus(\G(G))=1$ if and only if $\Bc(G)$ has exactly one edge, i.e.\ $G$ has exactly one cyclic subgroup of order $6$.
	\end{proof}
	
	Note that the proofs of Theorems~\ref{thm:planar} and~\ref{thm:torus} use only Kuratowski's theorem, Euler's formula and the counting of Lemma~\ref{lem:starlow}.
	
	\subsection{Genus two}\label{ss:genus2}
	
	For genus two we need some lower bounds of the form $\genus\ge3$. The first one concerns $\G(\Z_{12})$, which is not covered by Lemma~\ref{lem:starlow}.
	
	\begin{lemma}\label{lem:Z12}
		$\genus(\G(\Z_{12}))\ge3$.
	\end{lemma}
	\begin{proof}
		The parts are $A=\{e,a,a',s,s'\}$ ($a,a'$ of order $4$, $s,s'$ of order $6$) and $B=\{t,g_1,\dots,g_4,u,u^2\}$ ($t$ the involution, $g_i$ the generators, $u$ of order $3$), with $N(t)=A$, $N(g_i)=\{a,a',s,s'\}=:Q$, $N(u)=N(u^2)=\{e,s,s'\}$; $V=12$, $E=27$ and the minimum degree is $3$. Suppose $\genus=2$. Then $F=13$ and the face lengths sum to $54$, so there are twelve $4$-faces and one $6$-face. Apply Lemma~\ref{lem:balance} with the roles of the parts exchanged ($M$ on $A$; the dummy $z$ of the hexagon is a triangle on three distinct vertices of $A$, since a boundary walk does not backtrack at vertices of degree at least $2$). The faces of $M$ are: a pentagon $t$ through all of $A$; triangles $u,u^2$ on $\{e,s,s'\}$; $4$-cycles $g_1,\dots,g_4$ on $Q$; and the triangle $z$.
		
		If $u,u^2$ are oppositely oriented, they cancel on every pair; the two sides of $t$ at $e$ (say $x\to e\to y$) can only be balanced by $z$, so $z$ is the triangle $e\to x\to y\to e$; then $z$'s side $x\to y$ together with the remaining path $y\to p\to q\to x$ of $t$ is a directed $4$-cycle on $Q$. If $u,u^2$ have the same orientation, they contribute twice on each of the pairs at $e$, and only $t$ and $z$ can balance this; this forces $t$ to contain $s\to e\to s'$ and $z=(e,s',s)$, and again the net side on $\{s,s'\}$ together with the rest of $t$ is a directed $4$-cycle on $Q$. In both cases five directed $4$-cycles on the four vertices of $Q$ sum to zero on every pair. The $4$-cycles on four vertices are of three types, type $j$ avoiding the $j$-th perfect matching; a pair of matching $j$ is used by the cycles not of type $j$, so the number $n_j$ of cycles of type $j$ is odd for every $j$. For the matching $\{12,34\}$, types $2$ and $3$ (cycles $1234$ and $1243$) traverse $12$ in the same direction but $34$ in opposite directions; with $b,c$ the signed sums of the type-$2$ and type-$3$ cycles we get $b+c=0=b-c$, contradicting the oddness of $b$ and $c$.
	\end{proof}
	
	\begin{lemma}\label{lem:three}
		Let $K$ be a component of $\Bc(G)$ with an edge. Then $\genus(\G_K)\ge3$ in each of the following cases:
		\begin{enumerate}[label=(\roman*)]
			\item $K$ contains a subgroup of order $3$ of degree at least $5$;
			\item $K$ contains an involution of degree at least $7$;
			\item $K$ is not a star.
		\end{enumerate}
	\end{lemma}
	\begin{proof}
		(i) By Lemmas~\ref{lem:Bdeg} and~\ref{lem:no5}, such a subgroup $U$ has degree $m\ge7$. The vertex $e$, the two generators of $U$, seven involutions centralizing $U$ and the generators of the seven corresponding cyclic subgroups of order $6$ span a copy of $S_U(7)$ in $\G_K$, and $\genus(S_U(7))\ge\lceil29/12\rceil=3$ by Lemma~\ref{lem:starlow}(i).
		
		(ii) In the same way $\G_K\supseteq S_t(7)$, and $\genus(S_t(7))\ge\lceil7/3\rceil=3$ by Lemma~\ref{lem:starlow}(ii).
		
		(iii) By Lemma~\ref{lem:shape}(c), $K$ has an involution $t$ of degree at least $4$ adjacent to a subgroup $U_1$ of degree at least $3$. By (i), (ii) and Lemma~\ref{lem:Bdeg} we may assume that $t$ has degree exactly $4$ and $U_1$ degree exactly $3$. Let $L=C_G(t)$, let $U_1,\dots,U_4$ be the subgroups of order $3$ of $L$, and let $d_i$ be the degree of $U_i$ in $\Bc(G)$.
		
		\emph{\textbf{Step 1:} $d_i\ge3$ for all $i$.} If $y\in L$ has order $9$, then $ty$ would have order $18\notin\Pi$; so a Sylow $3$-subgroup $T$ of $L$ has exponent $3$. A group of exponent $3$ and order at least $27$ has at least $13$ subgroups of order $3$, and two distinct subgroups of order $9$ and exponent $3$ contain at least $4+4-1=7$ subgroups of order $3$; since $L$ has exactly four, either $|T|=3$, or $T\cong\Z_3^2$ is the unique Sylow $3$-subgroup of $L$ and contains $U_1,\dots,U_4$. In the first case $U_1,\dots,U_4$ are the Sylow $3$-subgroups of $L$, so $U_i=lU_1l^{-1}$ for some $l\in L$; conjugation by $l$ maps $C_G(U_1)$ onto $C_G(U_i)$, and $d_i=d_1=3$. In the second case $T$ is abelian, so $T\le C_G(U_1)$, and $T$ acts by conjugation on the three involutions of $C_G(U_1)$, fixing $t$. The orbits of a $3$-group have size $1$ or $3$, so $T$ fixes all three, i.e.\ the three involutions of $C_G(U_1)$ centralize $T\supseteq U_i$, and $d_i\ge3$.
		
		\emph{\textbf{Step 2:} counting.} Let $I$ be the set of involutions adjacent in $K$ to at least one $U_i$, and let $\varepsilon=\sum_id_i\ge12$ be the number of edges of $K$ at $U_1,\dots,U_4$. Since $t$ is adjacent to every $U_i$, $|I|\le1+\sum_i(d_i-1)=\varepsilon-3$. Let $Y\subseteq\G_K$ be the subgraph consisting of $e$, the elements of $I$, the eight generators of $U_1,\dots,U_4$, and, for each of the $\varepsilon$ edges $sU_i$, the two generators of $\langle s\rangle\times U_i$; its edges join $e$ to the elements of $I$ and to the eight generators, and each element of order $6$ to its cube $s$ and to the two generators of $U_i$ (indices $3$ and $2$). Then $Y$ is connected and bipartite, with $V=9+|I|+2\varepsilon$ and $E=8+|I|+6\varepsilon$. By Fact~\ref{fact:graph}(v),
		\[ \genus(\G_K)\ \ge\ \genus(Y)\ \ge\ \Big\lceil\frac{E-2V+4}{4}\Big\rceil=\Big\lceil\frac{2\varepsilon-|I|-6}{4}\Big\rceil\ \ge\ \Big\lceil\frac{\varepsilon-3}{4}\Big\rceil\ \ge\ \Big\lceil\frac94\Big\rceil=3. \qedhere \]
	\end{proof}
	
	\begin{theorem}[Genus two]\label{thm:genus2}
		Let $G$ be a finite group. Then $\genus(\G(G))=2$ if and only if every element order of $G$ lies in $\Pi$ and one of the following holds:
		\begin{enumerate}[label=(\alph*)]
			\item $G$ has exactly one subgroup $U$ of order $3$ that is centralized by an involution, and $C_G(U)$ contains exactly three involutions; or
			\item $G$ has exactly one involution $t$ that centralizes an element of order $3$, and $C_G(t)$ contains exactly eight elements of order $3$ (equivalently, exactly four subgroups of order $3$).
		\end{enumerate}
		Equivalently: all element orders lie in $\Pi$, and $\Bc(G)$ has exactly one component with an edge, which is a star centred at a subgroup of order $3$ with three involutions, or a star centred at an involution with four subgroups of order $3$. Elements of order $4$, $8$ and $9$ are allowed in both cases.
	\end{theorem}
	\begin{proof}
		\emph{Sufficient condition.} In case (a), every edge of $\Bc(G)$ contains $U$, and $U$ has degree $3$; so $\Bc(G)$ has exactly one component with an edge, a star centred at $U$ with three involutions. In case (b), every edge contains $t$, which has degree $4$, and the component is a star centred at $t$ with four subgroups of order $3$. By Lemma~\ref{lem:shape}(a),(b), $\G_K$ is $S_U(3)$, resp.\ $S_t(4)$, together with vertices of order $9$ adjacent exactly to the two generators of a subgroup $U_i$; by Lemma~\ref{lem:starup}, $\genus(\G_K)=2$. By Lemma~\ref{lem:decomp}, $\genus(\G(G))=2$.
		
		\emph{Necessary condition.} Suppose $\genus(\G(G))=2$. By Lemma~\ref{lem:orders} all element orders lie in $\Pi$, and by Lemma~\ref{lem:decomp} and the proof of Theorem~\ref{thm:torus}, either $\Bc(G)$ has exactly two components with an edge, each a single edge, or exactly one such component $K$, with $\genus(\G_K)=2$.
		
		The first case is impossible. There $G$ has exactly two cyclic subgroups $C_1=\langle t_1\rangle\times U_1$ and $C_2=\langle t_2\rangle\times U_2$ of order $6$, lying in different components, so $U_1\ne U_2$. Conjugation permutes $\{C_1,C_2\}$, so $N=N_G(C_1)$ has index at most $2$ in $G$; a subgroup of index at most $2$ is normal and contains all elements of odd order, so $U_2\le N$. The subgroup $\langle t_1\rangle$ is characteristic in $C_1\trianglelefteq N$, hence normal in $N$, and a normal subgroup of order $2$ is central. Thus $t_1$ centralizes $U_2$, i.e.\ $t_1\sim U_2$ in $\Bc(G)$, and $U_1$, $U_2$ lie in the same component, a contradiction.
		
		In the second case $K$ has at least two edges, since a single edge gives genus $1$. By Lemma~\ref{lem:three}, $K$ is a star, and its centre is a subgroup of order $3$ of degree less than $5$, hence of degree $3$ (Lemma~\ref{lem:Bdeg}), or an involution of degree less than $7$, hence of degree $4$. Hence $K$ is a star centred at a subgroup $U$ of order $3$ with exactly three involutions, or a star centred at an involution $t$ with exactly four subgroups of order $3$. In the first case every involution centralizing an element of order $3$ lies in $K$, hence centralizes only $U$; so $U$ is the only subgroup of order $3$ centralized by an involution, and $C_G(U)$ contains exactly three involutions, which is (a). In the second case $t$ is the only involution centralizing an element of order $3$, and the four subgroups of order $3$ centralized by $t$ give exactly eight elements of order $3$ in $C_G(t)$, which is (b).
	\end{proof}
	
	\section*{Acknowledgments}
	
	The research of Pallabi Manna is supported by Anusandhan National Research Foundation
	(ANRF)-India under the National Postdoctoral Fellowship scheme (File No.\ PDF/2025/006751). This 
	author acknowledges the excellent research facilities provided by the Department of Mathematics.
	The	second author acknowledges excellent research facilities provided by SRM University-AP. 
	
	\section*{Statements and Declarations}
	
	\noindent\textbf{Competing Interests.} The authors have no competing interests to declare that are
	relevant to the content of this article.
	
	\section*{Data Availability}
	
	Data sharing is not applicable to this article, as no datasets were used or generated during this
	study.
	


\begin{thebibliography}{99}
		
		\bibitem{AACNS} G.~Aalipour, S.~Akbari, P.~J.~Cameron, R.~Nikandish and F.~Shaveisi, On the structure of the power graph and the enhanced power graph of a group, \emph{Electron. J. Combin.} \textbf{24}(3) (2017), Article No. P3.16.
		
		\bibitem{AAM} A.~Abdollahi, S.~Akbari and H.~R.~Maimani, Non-commuting graph of a group, \emph{J. Algebra} \textbf{298} (2006), 468--492.
		
		\bibitem{AAHRS} S.~Akbari, A.~Ashtab, F.~Heydari, M.~Rezaee and F.~Sherafati, The prime index graph of a group, preprint, arXiv:1508.01133 (2015).
		
		\bibitem{AR} T.~Anitha and R.~Rajkumar, Prime index elements graph of a group, \emph{Comm. Algebra} \textbf{51}(3) (2023), 1250--1263.
		
		\bibitem{BHKY} J.~Battle, F.~Harary, Y.~Kodama and J.~W.~T.~Youngs, Additivity of the genus of a graph, \emph{Bull. Amer. Math. Soc.} \textbf{68} (1962), 565--568.
		
		\bibitem{BF} R.~Brauer and K.~A.~Fowler, On groups of even order, \emph{Ann. Math.} \textbf{62} (1955), 565--583.
		
		\bibitem{Cameron} P.~J.~Cameron, Graphs defined on groups, \emph{Int. J. Group Theory} \textbf{11}(2) (2022), 53--107.
		
		\bibitem{CGS} I.~Chakrabarty, S.~Ghosh and M.~K.~Sen, Undirected power graphs of semigroups, \emph{Semigroup Forum} \textbf{78} (2009), 410--426.
		
		\bibitem{CEJO} G.~Chartrand, L.~Eroh, M.~A.~Johnson and O.~R.~Oellermann, Resolvability in graphs and the metric dimension of a graph, \emph{Discrete Appl. Math.} \textbf{105} (2000), 99--113.
		
		\bibitem{Gor} D.~Gorenstein, \emph{Finite Groups}, Harper \& Row, New York, 1968.
		
		\bibitem{Hall} M.~Hall, \emph{The Theory of Groups}, Macmillan, New York, 1959.
		
		\bibitem{KQ} A.~V.~Kelarev and S.~J.~Quinn, Directed graphs and combinatorial properties of semigroups, \emph{J. Algebra} \textbf{251} (2002), 16--26.
		
		\bibitem{KSCT} A.~Kumar, L.~Selvaganesh, P.~J.~Cameron and T.~Tamizh Chelvam, Recent developments on the power graph of finite groups -- a survey, \emph{AKCE Int. J. Graphs Comb.} \textbf{18}(2) (2021), 65--94.
		
		\bibitem{OPF} O.~R.~Oellermann and J.~Peters-Fransen, The strong metric dimension of graphs and digraphs, \emph{Discrete Appl. Math.} \textbf{155} (2007), 356--364.
		
		\bibitem{RA17} R.~Rajkumar and T.~Anitha, Reduced power graph of a group, \emph{Electron. Notes Discrete Math.} \textbf{63} (2017), 69--76.
		
		\bibitem{Ringel} G.~Ringel, Das Geschlecht des vollst\"andigen paaren Graphen, \emph{Abh. Math. Sem. Univ. Hamburg} \textbf{28} (1965), 139--150.
		
		\bibitem{Williams} J.~S.~Williams, Prime graph components of finite groups, \emph{J. Algebra} \textbf{69} (1981), 487--513.
		
		\bibitem{Youngs} J.~W.~T.~Youngs, Minimal imbeddings and the genus of a graph, \emph{J. Math. Mech.} \textbf{12} (1963), 303--315.
		
	\end{thebibliography}
\end{document}